\documentclass{amsart}

\usepackage{amssymb,latexsym}
\usepackage{amsmath}
\usepackage{pb-diagram}
\usepackage{pstricks}
\usepackage{pst-tree}
\usepackage{mathrsfs}
\usepackage[all]{xy}

\newcommand{\id}{\operatorname{Id}}

\newcommand{\coder}{\operatorname{Coder}}

 \newtheorem{thm}{Theorem}
 \newtheorem{lem}[thm]{Lemma}

 \theoremstyle{definition}
 \newtheorem{defn}[thm]{Definition}
 \newtheorem*{ack}{Acknowledgments}

 \newtheorem{expl}[thm]{Examples}

\begin{document}

\title[Hochschild Cohomology is a BV Algebra]
{The BV Algebra on Hochschild Cohomology Induced by Infinity Inner
Products}
\author{Thomas Tradler}

\address{Thomas Tradler, College of Technology of the City University
of New York, Department of Mathematics, 300 Jay Street, Brooklyn, NY
11201, USA}

\email{ttradler@citytech.cuny.edu}

\begin{abstract}
We define a BV-structure on the Hochschild-cohomology of a unital, associative  algebra $A$ with a symmetric, invariant and non-degenerate inner product. The induced Gerstenhaber algebra is the one described in Gerstenhaber's original paper on Hochschild-cohomology. We also prove the corresponding theorem in the homotopy case, namely we define the BV-structure on the Hochschild-cohomology of a unital $A_\infty$-algebra with a symmetric and non-degenerate $\infty$-inner product.
\end{abstract}

\maketitle

\tableofcontents

\section{Introduction}

M. Gerstenhaber showed in \cite{G}, that for any associative algebra $A$, the Hochschild-cohomology $H^\bullet(A,A)$ has a Gerstenhaber-structure. More precisely, $H^\bullet(A,A)$ has a Gerstenhaber-bracket $[.,.]$ which is a Lie-bracket of degree $-1$, and a commutative, associative cup-product $\smile$, such that $[.,.]$ is a graded derivation of the cup-product in each variable. It is an important case, when the Gerstenhaber-structure comes from a Batalin-Vilkovisky (BV) algebra. This means that there is a degree $-1$ operator $\Delta$ on $H^\bullet(A,A)$, such that $\Delta$ squares to zero, and so that the deviation of $\Delta$ from being a derivation is the bracket.

Our goal in this paper is to define such a $\Delta$-operator on $H^\bullet(A,A)$, under the additional condition that $A$ has a unit $1$ and a symmetric, invariant and non-degenerate inner product. We will in fact produce two versions of this theorem; one for associative algebras, and one for homotopy associative algebras. To be precise about our first version, we now state the theorem, that is proved in section \ref{Assoc-BV}.
\begin{thm}\label{strict-thm}
Let $A$ be a finite dimensional, ungraded, unital, associative algebra with symmetric, invariant and non-degenrate inner product $<.,.>A\otimes A\to R$. For $f\in C^n(A,A)$, let $\Delta f\in C^{n-1}(A,A)$ be given by the equation
\begin{equation*}
<\Delta f (a_{1},\dots, a_{n-1}),a_{n}>
=\sum_{i=1}^n (-1)^{i(n-1)} <f(a_{i}, \dots, a_{n-1}, a_{n}, a_{1},\dots, a_{i-1}),1>.
\end{equation*}
Then, $\Delta$ is a chain map, such that the induced operation on Hochschild-cohomology $H^\bullet(A,A)$ squares to zero. Furthermore, for $\alpha\in H^n(A,A)$, $\beta\in H^m(A,A)$, the following identity holds,
\begin{equation*}
 [\alpha,\beta]= -(-1)^{(n-1)m} \cdot \Big( \Delta(\alpha\smile \beta)
  - \Delta(\alpha)\smile \beta  -(-1)^{n} \alpha\smile\Delta(\beta)\Big).
\end{equation*}
Here, $[.,.]$ and $\smile$ denote the Gerstenhaber-bracket and the cup-product defined in \cite{G}. 
\end{thm}
Is it important to remark, that the $\Delta$-operator in the above theorem is just Connes' B-operator on the Hochschild chains $C_\bullet(A,A)$, see \cite{C} or \cite[(2.1.9.1)]{L}, dualized and transferred to the Hochschild cochains $C^\bullet(A,A)$ via the inner product:
\begin{equation*}
\xymatrix{  C^\bullet(A,A)\ar[d] &&\text{Gerstenhaber-structure } [.,.], \smile \\
C^\bullet(A,A^*) \ar[u]_{\text{isomorphism induced by } <.,.>} && \text{dual of Connes' B-operator}}
\end{equation*}

Although this version of the theorem is already useful in certain situations, it is dissatisfactory in more general cases. In the second version of our theorem, we therefore remove the assumptions of finite dimensionality, we include grading, and use infinity versions of the concepts in theorem \ref{strict-thm}. More precisely, we start with an $A_\infty$-algebra $A$, and we use the concept of $\infty$-inner products, which was developed by the author in \cite{T}. Following Stasheff's approach to $A_\infty$-algebras and its Hochschild-cochain complex, see \cite{S} and \cite{S2}, using the tensor coalgebra and its coderivation complex, we define $A_\infty$-bimodules and the Hochschild-cochain complex with values in an $A_\infty$-bimodule, via the tensor bi-comodule, and coderivations associated to it. For an $A_\infty$-algebra $A$, $\infty$-inner products are $A_\infty$-bimodule maps between the $A_\infty$-bimodules $A$ and its dual $A^*$, see section \ref{infty-ip-section}.

It is known, that Hochschild-cohomology of an $A_\infty$-algebra still has a Gerstenhaber-structure, see \cite{GJ2}. In section \ref{infinity-case-section}, we generalize theorem \ref{strict-thm} to the infinity case as follows.
\begin{thm}\label{2x44}
Let $(A,D)$ be an A$_\infty$-algebra with strict unit $1$ and let $F$ be a symmetric, non-degenerate $\infty$-inner product. Then, $\Delta$ and the generalized cup-product $M$ induce a BV-algebra on the Hochschild-cohomology $H^\bullet(A,A)$. In fact, for $\alpha, \beta\in H^\bullet(A,A)$, we have that
$$ [\alpha,\beta]= \Delta(M(\alpha,\beta))-M(\Delta \alpha, \beta) -(-1)^{||\alpha||} M( \alpha,\Delta \beta), $$
where $[.,.]$ denotes, up to sign, the Gerstenhaber-bracket from \cite{GJ2}.
\end{thm}
We remark, that the difference of signs in theorems \ref{strict-thm} and \ref{2x44} stem from the difference of signs of the bracket in \cite{G} and the bracket used in theorem \ref{2x44}.

The ideas for both theorem \ref{strict-thm} and \ref{2x44} are taken from the corresponding proof in M. Chas and D. Sullivan's proof of the BV-algebra in string topology \cite{CS}. In fact, there is a strong analogy between \cite[lemma 4.6]{CS} and \cite[theorem 8.5]{G}. In this paper, we transfered the corresponding BV-proof \cite[lemma 5.2]{CS} to the Hochschild-cochain level.

It is an interesting question to determine the exact relationship between the BV-algebra given here, and the one from string topology \cite{CS}. Using the idea of local constructions of infinity structures from Lawrence and Sullivan in \cite{S3}, M. Zeinalian and the author constructed in \cite{TZ} a symmetric, non-degenerate $\infty$-inner product on the cochains $A=C^\bullet(X)$ of a Poincar\'e duality space $X$. If $X$ is simply connected, then it is known that the Hochschild-cohomology $H^\bullet(A,A)$ is isomorphic to the homology of the free loop space $H_\bullet(LX)$ of $X$, see \cite{J}. Thus, we can ask, if the BV-algebra on $H^\bullet(A,A)$ from theorem \ref{2x44}, and the one on $H_\bullet(LX)$ from string topology coincide under this isomorphism. We conjecture, that this is indeed the case, although we do not address this question in this paper.

\begin{ack}
I would like to thank Dennis Sullivan and Jim Stasheff for many helpful comments on this subject.
\end{ack}

\section{The associative case}\label{Assoc-BV}

In this section, we prove theorem \ref{strict-thm}. For this, we show that the BV-relations from theorem \ref{strict-thm} are satisfied on the Hochschild-cohomology of a unital, associative algebra with symmetric, invariant and non-degenerate inner product.

Let $R$ be a commutative ring such that the assumptions in \cite{G} are satisfied. Let $A$ be an (ungraded) finite dimensional associative algebra over $R$ with unit $1\in A$, and denote by $<.,.>:A\otimes A\to R$ a non-degenerate inner product such that $<a_1,a_2>=<a_2,a_1>$ and $<a_1\cdot a_2,a_3>=<a_1,a_2\cdot a_3>$ for all $a_1, a_2, a_3\in A$. The assumptions imply that the inner product induces a bimodule isomorphism between $A$ and its linear dual $A^*=Hom(A,R)$, given by $F:A\to A^*, a\mapsto F(a):=<a,.>$.

For any bimodule $M$ over $A$, denote by $C^n(A,M):=Hom(A^{ \otimes n},M)$, and $C^\bullet:=\bigoplus_{n\geq 0}C^n(A,A)$ the Hochschild-cochains of $A$ with values in $M$. For $f\in C^n(A,A)$, the differential $\delta f\in C^{n+1}(A,M)$ is given by
\begin{eqnarray*}
\delta f(a_1,\dots,a_{n+1})&:=& a_1. f(a_2,\dots,a_{n+1})\\
&& +\sum_{i=1}^n (-1)^i f(a_1,\dots,a_i\cdot a_{i+1},\dots, a_{n+1})\\
&& +(-1)^{n+1} f(a_1,\dots,a_n). a_{n+1}.
\end{eqnarray*}
It is $\delta^2=0$, and its homology $H^\bullet(A,M)$ is called Hochschild-cohomology of $A$ with values in $M$. We are mainly interested in the bimodules $M=A$ and $M=A^*$. Note that the map $F$ induces a chain map $F_\sharp:C^\bullet(A,A)\to C^\bullet(A,A^*), f\mapsto F_\sharp(f):=F\circ f$.

The cup-product for $f\in C^n(A,A)$ and $g\in C^m(A,A)$ is given by $f\smile g\in C^{n+m}(A,A)$, as $ f\smile g(a_1,\dots,a_{n+m}):=f(a_1,\dots, a_n)\cdot g(a_{n+1},\dots, a_{n+m})$.

\begin{defn} [$\Delta$-operator]
Let $n\geq 1$ and $f\in C^n(A,A)$. For $i\in\{1,\dots,n\}$ define $\Delta_{i} f\in  C^{n-1}(A,A)$ by the equation
\begin{equation*}
<\Delta_{i} f (a_{1},\dots, a_{n-1}),a_{n}>
= <f(a_{i}, \dots, a_{n-1}, a_{n}, a_{1},\dots, a_{i-1}),1>.
\end{equation*}
With this, define $\Delta:=\sum_{i=1}^{n} (-1)^{i(n-1)} \Delta_{i}:C^{n}(A,A)\to C^{n-1}(A,A)$. 

We also need the maps $\Delta^1$ and $\Delta^2$, which for $f\in C^n(A,A)$ and $g\in C^m(A,A)$ are elements $\Delta^1(f\otimes g), \Delta^2(f\otimes g)\in  C^{n+m-1}(A,A)$,  
\begin{eqnarray*}
& \Delta^{1}(f\otimes g):= \sum_{i=1}^{m} (-1)^{i(n+m-1)} \Delta_{i}(f\smile g),\\
 &  \Delta^{2}(f\otimes g):= \sum_{i=m+1}^{n+m} (-1)^{i(n+m-1)} \Delta_{i}(f\smile g). 
\end{eqnarray*}
We have $\Delta^{1} (f\otimes g)+ \Delta^{2} (f\otimes g)= \Delta(f\smile g)$. The notation for $\Delta^{1}$ and $\Delta^{2}$ is taken from \cite[chapter 5]{CS}.
\end{defn}

The following lemma follows by a straightforward calculation.
\begin{lem}\label{hoch-lemma}
In the above notation, we have that
\begin{enumerate}
\item $\Delta$ is a chain map,
\item\label{Delta1=Delta2} $\Delta^{1}(f\otimes g)=(-1)^{nm}\Delta^{2}(g\otimes f) $ for all $f\in C^{n}(A,A), g\in C^{m}(A,A)$.
\end{enumerate}
\end{lem}

Recall from \cite{G}, that for $j\in\{1,\dots,n\}$ and $f\in C^n(A,A), g\in C^m(A,A)$, $f\circ_{j} g\in C^{n+m-1}(A,A)$ is defined by $$ f\circ_{j} g (a_{1},\dots, a_{n+m-1}):=f(a_{1},\dots, a_{j-1} , g(a_{j},\dots, a_{j+m-1}),a_{j+m},\dots,a_{n+m-1}). $$
From this, one defines for $f\in C^n(A,A), g\in C^m(A,A)$ the pre-Lie product $f\circ g\in C^{n+m-1}(A,A)$ and the Gerstenhaber-bracket $[f,g]\in C^{n+m-1}(A,A)$, 
\begin{eqnarray*}
 f\circ g &:=& \sum_{j=1}^n (-1)^{(j-1)(m-1)} f\circ_j g,\\
\,[\, f , g\,] &:=& f\circ g -(-1)^{(n-1)(m-1)} g\circ f.
\end{eqnarray*}

\begin{lem}\label{homotopy-H}
For closed elements $f\in C^{n}(A,A)$ and $g\in C^{m}(A,A)$, $\delta f=\delta g=0$, we define $H\in C^{n+m-2}(A,A)$ by
$$ H:=\sum_{i\geq 1, j\geq 1, i+j\leq n}  (-1)^{(j-1)(m-1)+i(n+m)+1}
\Delta_{i} (f\circ_{j}g). $$ Then $H$ satisfies $ \delta H= f\circ g-(-1)^{(n-1)m}  \Delta(f)\smile g +(-1)^{(n-1)m} \Delta^{2} (f\otimes g)$.
\end{lem}

Before proving lemma \ref{homotopy-H}, we first show how this implies theorem \ref{strict-thm}.
\begin{proof}[Proof of theorem \ref{strict-thm}] 
We need to show that the Gerstenhaber-bracket $[.,.]$ is the deviation of $\Delta$ from being a derivation of $\smile$, 
\begin{eqnarray*}
 [f,g]&=&f\circ g-(-1)^{(n-1)(m-1)} g\circ f\\
      &\cong& (-1)^{(n-1)m} \Delta(f)\smile g - (-1)^{(n-1)m} \Delta^{2} (f\otimes g)\\
&&              -(-1)^{m-1} \Delta(g)\smile f +(-1)^{m-1}\Delta^{2} (g\otimes f)\\
      &\cong & -(-1)^{(n-1)m} \Delta^{2} (f\otimes g) -(-1)^{(n-1)m} \Delta^{1} (f\otimes g)\\
&&             +(-1)^{(n-1)m}\Delta(f)\smile g -(-1)^{(n-1)(m-1)}f\smile\Delta(g)\\
      &=& -(-1)^{(n-1)m} \cdot \big(\Delta(f\smile g) -\Delta(f)\smile g  -(-1)^{n} f\smile\Delta(g)\big),
\end{eqnarray*}
where we have used lemma \ref{hoch-lemma}\eqref{Delta1=Delta2} and the fact that the cup-product is graded commutative on Hochschild-cohomology. Furthermore it is well-known, that on homology the square of $\Delta$ vanishes. This can be seen, just as for Connes' B-operator, by considering the normalized Hochschild-cochain complex, which is a quasi-isomorphic subcomplex of $C^\bullet(A,A)$, where $(\Delta)^2=0$. 
\end{proof}

\begin{proof}[Proof of lemma \ref{homotopy-H}] 
Let $a_{1}, \dots,a_{n+m}\in A$ be fixed. For $i\leq j$, we will use the notation from \cite{G}, and denote by $a_{i,j}=(a_{i},\dots, a_{j})\in A^{\otimes j-i+1}$. For $i\geq 1$, $j\geq 1$, $i+j\leq n$, we calculate $\delta(\Delta_i(f\circ_j g))$ as follows.
\begin{eqnarray*}
&&<\delta(\Delta_i(f \circ_j g))(a_{1,n+m-1}),a_{n+m}>\\
&=&<a_1\cdot \Delta_i(f \circ_j g)(a_{2,n+m-1}),a_{n+m}>\\
&& +\sum_{\lambda=1}^{n+m-2} (-1)^\lambda <\Delta_i(f\circ_j g)(a_{1,\lambda-1},(a_\lambda\cdot a_{\lambda+1}),a_{\lambda+2,n+m-1}),a_{n+m}>\\
&&+ (-1)^{n+m-1} <\Delta_i(f\circ_j g)(a_{1,n+m-2})\cdot a_{n+m-1},a_{n+m}>\\
\end{eqnarray*}
\begin{eqnarray*}
&=& <f(a_{i+1,i+j-1},g(a_{i+j,i+j+m-1}),a_{i+j+m,n+m-1},(a_{n+m}\cdot a_1),a_{2,i}),1>\\
&&+\sum_{\lambda=1}^{i-1}(-1)^\lambda <f(a_{i+1,i+j-1},
g(a_{i+j,i+j+m-1}) , \\
&&\quad\quad \quad \quad \quad\quad\quad\quad\quad\quad
a_{i+j+m,n+m}, a_{1,\lambda-1} ,(a_{\lambda}\cdot a_{\lambda+1}),
a_{\lambda+2,i}),1> \\
&&+\sum_{\lambda=i}^{i+j-2} (-1)^{\lambda} <f(a_{i,\lambda-1},(a_{\lambda}\cdot a_{\lambda+1}), a_{\lambda+2,i+j-1}, g(a_{i+j,i+j+m-1}) ,\\
&&\quad\quad\quad\quad\quad\quad\quad\quad \quad \quad
 a_{i+j+m,n+m} , a_{1,i-1}),1>\\
&&+\sum_{\lambda=i+j-1}^{i+j+m-2} (-1)^\lambda <f(a_{i,i+j-2},g(a_{i+j-1,\lambda-1},(a_\lambda\cdot a_{\lambda+1}),a_{\lambda+2,i+j+m-1}),\\
&&\quad\quad\quad\quad\quad\quad\quad\quad \quad \quad
a_{i+j+m,n+m},a_{1,i-1}),1> \\
& &+\sum_{\lambda=i+j+m-1}^{n+m-1} (-1)^\lambda <f(a_{i,i+j-2},
g(a_{i+j-1,i+j+m-2}) ,\\
&&\quad\quad\quad\quad\quad\quad\quad\quad \quad \quad
 a_{i+j+m-1,\lambda-1},(a_{\lambda}\cdot a_{\lambda+1})
, a_{\lambda+2,n+m} , a_{1,i-1}),1>,
\end{eqnarray*}
Using the fact that $\delta(g)=0$, we see that we can write this as
\begin{equation}\label{i,j}
<(\delta(\Delta_{i}(f\circ_{j} g))) (a_{1,n+m-1}),a_{n+m}>= h_{i,j}+h'_{i,j}+h''_{i,j},
\end{equation}
where we define
\begin{eqnarray*}
h_{i,j}&:=&(-1)^{i+1}<a_{i}\cdot f(a_{i+1,i+j-1},g(a_{i+j,i+j+m-1}), a_{i+j+m,n+m} , a_{1,i-1}),1>\\
& &+\sum_{\lambda=i}^{i+j-2} (-1)^{\lambda} <f(a_{i,\lambda-1},(a_{\lambda}\cdot a_{\lambda+1}), a_{\lambda+2,i+j-1}, g(a_{i+j,i+j+m-1}) ,\\
& & \quad\quad\quad\quad\quad\quad\quad\quad \quad
 a_{i+j+m,n+m} , a_{1,i-1}),1>\\
& &+(-1)^{i+j-1} <f(a_{i,i+j-2}, (a_{i+j-1}\cdot g(a_{i+j,i+j+m-1})), a_{i+j+m,n+m} , a_{1,i-1}),1>\\
h'_{i,j}&:=&(-1)^{i+j+m-2} <f(a_{i,i+j-2}, (g(a_{i+j-1,i+j+m-2})
\cdot a_{i+j+m-1}),\\
& & \quad\quad\quad\quad\quad\quad\quad\quad \quad
 a_{i+j+m,n+m} , a_{1,i-1}),1>\\
& &+\sum_{\lambda=i+j+m-1}^{n+m-1} (-1)^\lambda <f(a_{i,i+j-2},
g(a_{i+j-1,i+j+m-2}) , \\
& & \quad\quad\quad\quad\quad\quad\quad\quad \quad
a_{i+j+m-1,\lambda-1},(a_{\lambda}\cdot a_{\lambda+1})
, a_{\lambda+2,n+m} , a_{1,i-1}),1>,
\end{eqnarray*}
and
\begin{eqnarray*}
h''_{i,j}&:=&<f(a_{i+1,i+j-1}, g(a_{i+j,i+j+m-1}) ,
a_{i+j+m,n+m-1},(a_{n+m}\cdot a_{1}), a_{2,i}),1>\\
& &+\sum_{\lambda=1}^{i-1} (-1)^\lambda <f(a_{i+1,i+j-1},
g(a_{i+j,i+j+m-1}) , a_{i+j+m,n+m}, \\
& & \quad\quad\quad\quad\quad\quad\quad\quad \quad
a_{1,\lambda-1} ,(a_{\lambda}\cdot a_{\lambda+1}),
a_{\lambda+2,i}),1> \\
& &+(-1)^i <f(a_{i+1,i+j-1}, g(a_{i+j,i+j+m-1}) ,
a_{i+j+m,n+m}, a_{1,i-1})\cdot a_{i},1>.
\end{eqnarray*}
On the other hand, if $i\geq 2$, $j\geq 1$, $i+j\leq n-1$, we may use the fact that $\delta(f)=0$, to obtain
\begin{multline}\label{i,j+-1}
(-1)^{i+1} h_{i,j}+(-1)^{i+m} h'_{i,j+1}+(-1)^{i+n} h''_{i-1,j+1} \\ =<(\Delta_{i}((\delta f)
\circ_{j} g)) (a_{1,n+m-1}),a_{n+m}>=0.
\end{multline}
If we set furthermore for $i,j\in \{1,\dots,n\}$,
\begin{eqnarray*}
h_{i,0}&:=&(-1)^{i+1} <g(a_{i,i+m-1})\cdot f(a_{i+m,n+m}, a_{1,j-1}),
1>\\
h'_{i,n-i+1}&:=&(-1)^{n+m+1}<f(a_{i,n-1} ,(g(a_{n,n+m-1})\cdot a_{n+m})
, a_{1,i-1}),1>  \\
h''_{0,j}&:=&<f(a_{1,j-1}, g(a_{j,j+m-1}), a_{j+m,n+m-1}) \cdot a_{n+m},1>,
\end{eqnarray*}
then all three quantities $h_{i,j}$,
$h'_{i,j+1}$, $h''_{i-1,j+1}$ are defined for all $i, j$ with  $i\geq 1$, $j\geq 0$, $i+j\leq n$. Furthermore equation \eqref{i,j+-1} holds for each of these $i, j$, so that we obtain
\begin{multline*}
0=\sum_{i\geq 1, j\geq 0, i+j\leq n} (-1)^{(j-1)(m-1)+i(n+m-1)}\\
\cdot \big( (-1)^{i+1} h_{i,j} +(-1)^{i+m} h'_{i,j+1} +(-1)^{i+n} h''_{i-1,j+1} \big).
\end{multline*}
Rearranging the terms in this sum, and using equation \eqref{i,j}, gives
\begin{eqnarray*}
0&=&\sum_{i\geq 1, j\geq 1, i+j\leq n} (-1)^{(j-1)(m-1)+i(n+m)+1}\big(  h_{i,j} + h'_{i,j} + h''_{i,j} \big)\\
&& -\sum_{i=1}^{n} (-1)^{m-1+i(n+m)} h_{i,0} -
 \sum_{i=1}^{n} (-1)^{n(m+1)+i(n+1)} h'_{i,n-i+1} \\
 && - \sum_{j=1}^{n} (-1)^{(j-1)(m-1)} h''_{0,j} \\
 &=& <\delta(H)(a_{1,n+m-1}),a_{n+m}> \\
 &&- (-1)^{m(n+1)} <(\Delta^{2}(f\otimes g)) (a_{1,n+m-1}),a_{n+m}> \\
&& -(-1)^{m(n+1)+1} <(\Delta(f)\smile g)(a_{1,n+m-1}),a_{n+m}> \\
&&- <(f\circ g)(a_{1,n+m-1}),a_{n+m}>.
\end{eqnarray*}
This implies the claim of the lemma.
\end{proof}
In the case that $A$ is a graded associative algebra, we may obtain theorem \ref{strict-thm} in the same way as above, except for an appropriate change in signs. The next section includes graded associative algebras as a special case.

\section{The homotopy associative case}\label{infinity-case-section}

We now consider the more general situation from theorem \ref{2x44}. We recall the notion of $\infty$-inner products from \cite{T} in subsection \ref{infty-ip-section}, define the operations on the Hochschild-cochain complex in subsection \ref{def-operation-section}, and prove theorem \ref{2x44} in subsection \ref{proof-bv-section}. 

\subsection{Infinity-inner products}\label{infty-ip-section}
 
This section briefly recalls the definitions of
A$_\infty$-algebra, $\infty$-inner-products and the
Hochschild-cochain complex as defined in \cite{T}.

\begin{defn}[$A_\infty$-algebra]
 Let $A=\bigoplus_{j\in
\mathbb{Z}} A_{j}$ be a graded module over a given ground ring
$R$. Define its suspension $sA$ to be the graded module
$sA= \bigoplus_{j\in \mathbb{Z}} (sA)_{j}$ with $(sA)_{j}:=
A_{j-1}$. The suspension map $s:A\longrightarrow sA$, $a\mapsto
sa:=a$ is an isomorphism of degree +1.
We define $BA:=TsA=\bigoplus_{i\geq 0} sA^{\otimes i}$  to be the tensor coalgebra of $sA$, with the comultiplication
  $$
  \Delta:BA\longrightarrow BA\otimes BA,
  \quad
  \Delta(a_{1},
  \ldots  ,a_{n}):=\sum_{i=0}^{n} (a_{1},\ldots,a_{i})
  \otimes(a_{i+1},\ldots,a_{n}).
  $$
A coderivation on $BA$ is a map $D:BA\to BA$ such that
  $$
\begin{diagram}
  \node{BA}\arrow{e,t}{\Delta} \arrow{s,l}{D}
  \node{BA\otimes BA}\arrow{s,r}{D\otimes \id+\id\otimes D}\\
  \node{BA}\arrow{e,b}{\Delta} \node{BA\otimes BA}
\end{diagram}
$$
The space of coderivation on $BA$ is denoted by $Coder(BA)$.
By definition, an A$_\infty$ algebra structure on $A$ is a
coderivation $D\in \coder(BA)$ of degree $-1$ with $D^{2}=0$.

Every $A_\infty$-algebra $(A,D)$ is uniquely given by its components $D=\sum_{i\geq 0} D_i$, where $D_i:sA^{\otimes i}\to sA$ are lifted to coderivations, which, by abuse of notation, are still denoted by $D_i$, cf. \cite[section 2]{T}. We follow the usual convention and assume that the lowest level $D_0$ vanishes. Graphically, we represent $D$ by a filled circle, which takes many inputs and gives one output, here from left to right:
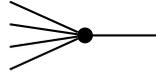
\[ \pstree[treemode=L, levelsep=1cm, treesep=0.3cm] {\Tp}{\pstree{\Tc*{3pt}}{\Tp \Tp \Tp \Tp}} \]
\end{defn}

\begin{defn}[$A_\infty$-bimodule]\label{A-infty-bimodule}
Now, let $(A,D)$ be an $A_\infty$-algebra, and let $M= \bigoplus_{j\in \mathbb{Z}}  M_j$ be a graded module. Define $B^{M}A:=R\oplus\bigoplus_{k\geq 0, l\geq 0} sA^{\otimes k} \otimes sM
  \otimes sA^{\otimes l}$ to be the tensor bi-comodule with the bi-comodule map
  $$
  \Delta^{M}:B^{M}A\longrightarrow(BA\otimes B^{M}A)\oplus (B^{M}A\otimes BA),$$
\begin{multline*}
  \Delta^{M}(a_{1},\ldots,a_{k},m,a_{k+1},\ldots,a_{k+l}):=
  \sum_{i=0}^{k}(a_{1},\ldots,a_{i}) \otimes(a_{i+1},\ldots,m,\ldots,a_{n}) \\
  +\sum_{i=k}^{k+l} (a_{1},\ldots,m,\ldots,a_{i})  \otimes(a_{i+1},\ldots,a_{k+l}).
\end{multline*}
Now, define a coderivation on $B^M A$ over $D$ to be a map $D^M:B^M A\rightarrow B^M A$
such that
$$
\begin{diagram}
  \node{B^M A}\arrow{e,t}{\Delta^{M}} \arrow{s,l}{D^M} \node{(BA\otimes
    B^M A)\oplus (B^M A\otimes BA)}\arrow{s,r}{(D\otimes \id+\id\otimes
    D^M)\oplus (D^M\otimes \id+\id\otimes D)}\\
  \node{B^M A}\arrow{e,b}{\Delta^{M}} \node{(BA\otimes B^M A)\oplus(B^M A\otimes BA)}
\end{diagram}
$$
The space of all coderivations on $B^M A$ over $D$ is denoted by
$\coder_{D}(B^M A)$. An A$_\infty$-bimodule structure on $M$ over $A$ is
defined to be a coderivation $D^{M}\in \coder_{D}(B^{M}A)$ of degree $-1$ with $(D^{M})^{2}=0$.

Finally, a coderivation from $BA$ to $B^M A$ is a map $f:BA\to B^M A$ such that
\[
\begin{diagram}
\node{BA}\arrow{e,t}{\Delta} \arrow{s,l}{f}
\node{BA\otimes BA}\arrow{s,r}{\id\otimes f+f\otimes \id}\\
\node{B^M A}\arrow{e,b}{\Delta^{M}} \node{(BA\otimes B^M A)\oplus (B^M A\otimes BA)}
\end{diagram}
\]
The space of all coderivations from $BA$ to $B^M A$ is denoted by
$C^\bullet(A,M)$, and is called the Hochschild-cochain complex of $A$ with values in $M$. $C^{\bullet}(A,M)$ has a differential $\delta^{M}(f) :=D^{M}\circ
f-(-1)^{|f|}f\circ D$ satisfying $(\delta^{M})^{2}=0$.
The homology $H^{\bullet}(A,M)$ is called the Hochschild-cohomology of $A$ with values in $M$.

Similar to before, we can uniquely decompose an element $f\in C^\bullet(A,M)$ in its components $f=\sum_{i\geq 0} f_i$, where $f_i:sA^{\otimes i}\to sM$ is lifted to a coderivation $f_i:BA\to B^M A$, see \cite[section 3]{T}. Here is the graphical representation for $f$, where the inputs on the left are from $sA$ and the output on the right lies in $sM$:
\[ \pstree[treemode=L, levelsep=1cm, treesep=0.3cm] {\Tp}{\pstree{\Tr{f}}{\Tp \Tp \Tp \Tp}} \]
\end{defn}

\begin{defn} [$\infty$-inner-product]\label{infty-defs}
Let $(A,D)$ be an A$_\infty$ algebra and let $(M,D^{M})$ and $(N,D^{N})$ be two A$_\infty$ bimodules over $A$. An $A_\infty$-bimodule map between $M$ and $N$ is a map $F:B^M A\to B^N A$ such that
  $$
\begin{diagram}
  \node{B^M A}\arrow{e,t}{\Delta^{M}} \arrow{s,l}{F} \node{(BA\otimes
    B^M A)\oplus (B^M A\otimes BA)}\arrow{s,r}{(\id\otimes F)\oplus
    (F\otimes \id)}\\
  \node{B^N A}\arrow{e,b}{\Delta^{N}} \node{(BA\otimes B^N A)\oplus(B^N A\otimes BA)}
\end{diagram}
$$
and  $ D^{N}\circ F=F\circ D^{M}$.
Every $A_\infty$-bimodule map
$F:B^{M}A\to B^{N}A$ induces a chain map on the Hochschild-cochain complexes
$F_{\sharp}:C^{\bullet}(A,M) \to C^{\bullet}(A,N)$ by $f\mapsto F\circ f$.

Let $(A,D)$ be an $A_\infty$ algebra. Then both $A$ and its linear dual $A^*=Hom(A,R)$ have induced $A_\infty$-bimodule structures $D^{A}\in \coder_{D}(B^{A}A)$ and $D^{A^{*}}\in \coder_{D} (B^{A^{*}}A)$, see \cite[Lemma 3.9]{T}. An $\infty$-inner-product is by definition an $A_\infty$-bimodule map $F:B^A A\to B^{A^*} A$ between $A$ and $A^{*}$ of degree $0$. The $\infty$-inner-product is called non-degenerate if there exists an $A_\infty$-bimodule map $G:B^{A^{*}}A\to B^{A}A$, such that $F_\sharp:C^\bullet(A,A) \to C^\bullet(A,A^*)$ and $G_\sharp:C^\bullet(A,A^*) \to C^\bullet(A,A)$ are quasi-inverse to each other, so that their induced maps on Hochschild-cohomology are inverses.

Again, one can show, that the $\infty$-inner product $F$ is given by components $F:\sum_{i,j\geq 0} F_{i,j}$, where the $F_{i,j}:sA^{\otimes i}\otimes sA\otimes sA^{\otimes j}\to sA^*$ are lifted according to the property in definition \ref{infty-defs}, see \cite[section 4]{T}. If we apply arguments to this map, $F_{i,j}(a_1,\dots,a_i,a_{i+1},a_{i+2},\dots,a_{i+j+1})(a_{i+j+2})\in R$, then we can represent the $\infty$-inner product graphically as in \cite[section 5]{T} by an open circle:
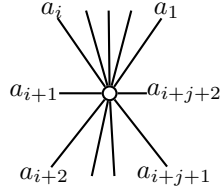
\[
  \pstree[treemode=R, levelsep=1cm, treesep=0.3cm]{\Tr*{a_{i+1}}}
    { \pstree[levelsep=0cm]{\Tc{3pt}}
  { \pstree[treemode=U, levelsep=0.8cm]{\Tn}
     {\Tr*{a_i} \Tp \Tp \Tp \Tr*{a_1}}
    \pstree[treemode=R, levelsep=1cm]{\Tn}
     {\Tr*{a_{i+j+2}} }
    \pstree[treemode=D, levelsep=0.8cm]{\Tn}
     {\Tr*{a_{i+2}} \Tp \Tp \Tr*{a_{i+j+1}}}
  }}
\]
Here we included the arguments $a_1,\dots, a_{i+j+2}$ to show their chosen position; the first $i$ arguments are placed on top, the special element $a_{i+1}$ on the left, the next $j$ arguments on the bottom and the output argument on the right. 

The $\infty$-inner product is called symmetric, if it is invariant under a $180^\circ$ rotation,
\begin{multline*}
\forall i,j:\quad F_{i,j}(a_1,\dots,a_i,a_{i+1},a_{i+2},\dots,a_{i+j+1})(a_{i+j+2})\\
= (-1)^{(\sum_{k=1}^{i+1}(|a_k|+1))(\sum_{k=i+2}^{i+j+2}(|a_k|+1))} F_{j,i}(a_{i+2},\dots,a_{i+j+1},a_{i+j+2},a_1,\dots,a_i)(a_{i+1}).
\end{multline*}
In our generalization of theorem \ref{strict-thm}, we use a symmetric, non-degenerate $\infty$-inner product in the above sense. We call such a structure $F$ an $\infty$-Poincar\'e-duality structure.
\end{defn}

\subsection{Definition of the operations}\label{def-operation-section}

We now describe the operations that are needed to obtain the BV-algebra on the  Hochschild-cohomology of an $A_\infty$-algebra. To this end, we first recall the Gerstenhaber-bracket and the cup-product from \cite{GJ2}, and the $\Delta$-operator from the previous section. These operations will only be up to sign the needed operations, and we denote them by $[.,.]'$, $M'$ and $\Delta'$. With the help of the $\infty$-inner product, we can rewrite these operations as maps of the form $C^\bullet(A,A)^{\otimes p}\to C^\bullet(A,A^*)$. Using a graphical representation of a large class of operations $C^\bullet(A,A)^{\otimes p}\to C^\bullet(A,A^*)$, which we call symbols, we may find symbols which represent the above operations up to sign. The operations of these symbols, toegther with their naturally given signs will be denoted by $[.,.]$, $M$ and $\Delta$, and will be used in subsection \ref{proof-bv-section} to prove theorem \ref{2x44}.

\begin{defn}[Gerstenhaber-bracket, cup-product]\label{gerst-cup-def}
Fix an $A_\infty$-algebra $(A,D)$, which in components is given by $D=\sum_{i\geq 1} D_i$, where $D_i:sA^{\otimes i}\to sA$. Let $f,g\in C^\bullet(A,A)$ be represented in components by $f=\sum_{n\geq 0} f_n$, and $g=\sum_{m\geq 0} g_m$ as in definition \ref{A-infty-bimodule}, where $f_n:sA^{\otimes n}\to sA$ and $g_m:sA^{\otimes m}\to sA$. For these components, we define the components of the Gerstenhaber-bracket $[f_n,g_m]':sA^{\otimes n+m-1}\to sA$ and the components of the generalized cup-product $M(f_n,g_m):\bigoplus_{k\geq 0} sA^{\otimes n+m+k}\to sA$ by
$$ [f_n,g_m]'=f_n\circ g_m -(-1)^{||f_n||\cdot ||g_m||} g_m\circ f_n, \text{ where }  $$
\begin{multline*}
f_n \circ g_m (a_1,\dots, a_{n+m-1})\\ 
=\sum_{i=1}^n (-1)^{||g_m||\cdot \sum_{l=1}^{i-1}(|a_l|+1)} f_n(a_1,\dots,g_m(a_i,\dots,a_{i+m-1}),\dots,a_{n+m-1}), 
\end{multline*}
\begin{multline*}
M'(f_n,g_m)(a_1,\dots,a_{n+m+k})=\sum_{ \substack{1\leq i \\ i+n\leq j \\ j\leq n+k } } (-1)^{ ||f_n||\cdot \sum_{l=1}^{i-1}(|a_l|+1) + ||g_m||\cdot (||f_n||+\sum_{l=1}^{j-1}(|a_l|+1)) } \\ \cdot D_{k-2}(a_1,\dots,f_n(a_i,\dots, a_{i+n-1}),\dots, g_m(a_j,\dots,a_{j+m-1}),\dots, a_{n+m+k}).
\end{multline*}
where $||f_n||$ and $||g_m||$ are the total degrees of the maps $f_n:sA^{\otimes n}\to sA$ and $g_m:sA^{\otimes m}\to sA$.
Both the bracket $[.,.]'$ and the product $M'$ are extended linearly to operations $C^\bullet(A,A)\otimes C^\bullet(A,A)\to C^\bullet(A,A)$.
\end{defn}
 The following is observed by E. Getzler and J. Jones in \cite{GJ2}.
\begin{thm}\label{2x38}
Let $(A,D)$ be an A$_\infty$-algebra. Then $[.,.]'$ and $M'$ induce a Gerstenhaber-structure on  Hochschild-cohomology $H^\bullet(A,A)$.
\end{thm}
Next, we consider the unit and the $\Delta'$-operator.
\begin{defn}[Unit, $\Delta'$-operator]
An element $1\in A_{0} \cong (sA)_{+1} \subset BA$ is called a strict unit of $A$, or simply a unit of $A$, if $D_n$ applied to any element of the form $(a_{1},..., a_{i-1},1,a_{i+1},...,a_{n})\in (sA)^{\otimes n}$ vanishes, except for the case $n=2$,
\begin{eqnarray*}
& D_n(a_{1},...,a_{i-1},1,a_{i+1},...,a_{n})=0 & \text{ for } n\neq 2, \\
& -D_2(1,a) = (-1)^{|a|+1} D_2(a,1)=a. &
\end{eqnarray*}
The sign is chosen so that the unshifted multiplication $m_2=s^{-1}\circ D_2\circ (s\otimes s)$ satisfies $m_2(1,a)=m_2(a,1)=a$, see \cite[proposition 2.4]{T}.

We also have the $\Delta'$-operator from section \ref{Assoc-BV}, which we write as an operation $\Delta':C^\bullet(A,A^*)\to C^\bullet(A,A^*)$, which is the dual of Connes' B-operator.  For $f\in C^\bullet(A,A^*)$ with components $f=\sum_{n\geq 0} f_n$, where $f_n:sA^{\otimes n}\to sA^*$, we define the component $\Delta' f_n :sA^{n-1}\to sA^*$ by
\begin{multline*}
\Delta' f_n(a_1,\dots, a_{n-1})(a_n) \\ = \sum_{i=1}^{n} (-1)^{\sum_{l=1}^{i-1}(|a_l|+1) \cdot \sum_{l=i}^{n}(|a_l|+1)} f_n(a_i,\dots,a_{n-1}, a_n,a_1,\dots, a_{i-1})(1).
\end{multline*}
It is straightforward to check that in the $A_\infty$ setting with strict unit $1$, $\Delta'$ still forms a chain map. Furthermore, the operation on Hochschild-cohomology induced by the square of $\Delta'$ vanishes, which can be seen by considering the normalized Hochschild-cochain subcomplex, where $(\Delta')^2=0$ holds identically.
\end{defn}

Now, fix a symmetric, non-degenrate $\infty$-inner product $F$. The non-degenracy of $F$ implies that $F_\sharp:C^\bullet(A,A)\to C^\bullet(A,A^*)$ induces an isomorphism of Hochschild-cohomologies $H^\bullet(A,A)\to H^\bullet(A,A^*)$. We use $F_\sharp:C^\bullet(A,A)\to C^\bullet(A,A^*)$ and its quasi-inverse $G_\sharp :C^\bullet(A,A^*)\to C^\bullet(A,A)$ from defintion \ref{infty-defs}, to write all necessary operations as maps $C^\bullet(A,A)^{\otimes p}\to C^\bullet(A,A^*)$. For example, we have the induced operations
\begin{eqnarray*}
[.,.]'&: & C^\bullet(A,A)^{\otimes 2}\stackrel{[.,.]'}{\longrightarrow} C^\bullet(A,A)\stackrel{F_\sharp}{\longrightarrow} C^\bullet(A,A^*), \\
M'&: & C^\bullet(A,A)^{\otimes 2}\stackrel{M'}{\longrightarrow} C^\bullet(A,A)\stackrel{F_\sharp}{\longrightarrow} C^\bullet(A,A^*), \\
\Delta' &: &  C^\bullet(A,A)\stackrel{F_\sharp}{\longrightarrow} C^\bullet(A,A^*)\stackrel{\Delta'}{\longrightarrow} C^\bullet(A,A^*),\\
\Delta'\circ M'&:&C^\bullet(A,A)^{\otimes 2}\stackrel{M'}{\longrightarrow} C^\bullet(A,A)\stackrel{F_\sharp}{\longrightarrow} C^\bullet(A,A^*)\stackrel{\Delta'}{\longrightarrow} C^\bullet(A,A^*), \\
M'\circ (\id\otimes \Delta') &: &  C^\bullet(A,A)^{\otimes 2}\stackrel{\id\otimes (G_\sharp\circ \Delta' \circ F_\sharp)}{\longrightarrow} C^\bullet(A,A)^{\otimes 2}\stackrel{M'}{\longrightarrow} C^\bullet(A,A)\stackrel{F_\sharp}{\longrightarrow} C^\bullet(A,A^*).
\end{eqnarray*}

We will now introduce a convenient graphical representation of a large class of operations $C^\bullet(A,A)^{\otimes p}\to C^\bullet(A,A^*)$.
\begin{defn}[Symbols]
A symbol is a finite, planar tree, where each vertex is labeled by one of the following labels:
\begin{itemize}
\item an open circle: \pstree[levelsep=0cm, treesep=0cm]{\Tc{3pt}}{\Tn}
\item a filled circle: \pstree[levelsep=0cm, treesep=0cm]{\Tc*{3pt}}{\Tn}
\item the unit $1$
\item a numbered circle: \pstree[levelsep=0cm, treesep=0cm]{\Tcircle{1}}{\Tn}, \pstree[levelsep=0cm, treesep=0cm]{\Tcircle{2}} {\Tn}, \pstree[levelsep=0cm, treesep=0cm]{\Tcircle{3}}
{\Tn}, $\dots$
\end{itemize}
We require the following restrictions for a symbol.
\begin{itemize}
\item 
There is exactly one vertex which is labeled by an open circle. This vertex has two distinguished edges attached, which in the plane will be drawn as horizontally incoming from the left and right. There are any amount of additional edges attached to this vertex, placed above or below the horizontal line.
\item
Each vertex labeled by a filled circle has at least three edges attached.
\item
Each vertex labeled by the unit has exactly one edge attached.
\item 
There is exactly one distinguished external edge in the symbol. This edge will be denoted by a fat edge:   \pstree[treemode=R, levelsep=1cm]{\Tp} {\Tp[edge={\ncline[linewidth=0.1cm]}]}
\end{itemize}

Every symbol $\sigma$ with numbered circles from $1$ to $p$ uniquely determines an operation $\rho_\sigma:C^\bullet(A,A)^{\otimes p}\to C^\bullet(A,A^*)$ in the following way. Let $f^1,\dots, f^p \in C^\bullet(A,A)$. We describe the element $\rho_\sigma(f^1,\dots, f^p)\in C^\bullet(A,A^*)$ by determining its value on arguments $\rho_\sigma(f^1,\dots, f^p)(a_1,\dots,a_n)(a_{n+1})\in R$. 

This value is obtained by cyclically attaching the arguments $a_1,\dots, a_{n+1}$ to vertices which are labeled by an open circle, a filled circle, or a numbered circle, such that the last argument $a_{n+1}$ lands on the distinguished external fat edge. Note, that each vertex labeled by a filled circle or a numbered circle has a unique path to the vertex labeled by the open circle. Thus there is one preferred edge at every filled and numbered circle, which lies on the path to the open circle. We call this edge the outgoing edge of this vertex, and all other edges attached to the vertex are considered as incoming to this vertex. We use this notion of incoming and outgoing edges at filled and numbered circles, to successively evaluate the arguments and the unit $1$  toward the open circle. More precisely, to every filled vertex, we evaluate the arguments and the incoming edges using the $A_\infty$-structure $D$, where $D:\bigoplus_{i\geq 1} sA^{\otimes i}\to sA$, and record the  answer of this evaluation at the outgoing edge. Similarly, to every vertex with the circle numbered by $r$, we evaluate using the Hochschild-cochain $f^r:\bigoplus_{i\geq 0} sA^{\otimes i}\to sA$. We finally end at the vertex labeled by the open circle, where we evaluate all remaining arguments and incoming edges using the $\infty$-inner product $F:\bigoplus_{j,i} sA^{\otimes i}\otimes sA\otimes sA^{\otimes j}\to sA^*$. Here, the horizontal edge from the left is identified with the special $(i+1)$th argument, and the horizontal edge from the right is identified with the argument plugged into $sA^*$. We thus obtain an overall number in $R$ for each possible way of attaching the arguments $a_1,\dots, a_{n+1}$.

With this, we define $\rho_\sigma(f^1,\dots, f^p)(a_1,\dots,a_n)(a_{n+1})$ to be the sum of all numbers in $R$, obtained by the above procedure from attaching the arguments $a_1,\dots, a_{n+1}$ in all possible ways, where each number is multiplied by the following sign factor. We use the sign $(+1)$ for the unique term which is obtained by applying all arguments $a_1,\dots, a_{n+1}$ to the position specified by the distinguished external fat edge. All other terms will be compared to this unique term, by looking at its linear expression using the $f^1,\dots,f^r$, $D$, and $F$, with arguments $a_1,\dots, a_{n+1}$ and units $1$ applied, such as
 $$ F_{i,j}(\dots,D(\dots, f^k(\dots)),\dots,f^m(\dots,f^n(\dots),\dots), \dots) \in R. $$
We then apply the usual Koszul sign rule of a factor of $(-1)^{|x|\cdot |y|}$ for each variable $x$ that has moved across the variable $y$ when comparing it to the expansion of the above term with sign $(+1)$.
\end{defn}
The following examples show some symbols and their associated operations.
\begin{expl}\label{symbol-examples}
All examples except for the last one will be studied without signs.
\begin{enumerate}
\item
Denote by $\sigma$ the following symbol:
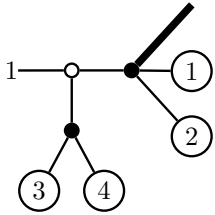
\[
   \pstree[treemode=R, levelsep=0.8cm, treesep=0.3cm,
     linewidth=0.9pt] {\Tr*{1}}{\pstree[levelsep=0cm]{\Tc{3pt}}
  { \Tn  \pstree[treemode=R, levelsep=0.8cm]{\Tn}
       {\pstree{\Tc*{3pt}} {\Tp[edge={\ncline[linewidth=0.1cm]}] \Tn \Tcircle{1}\Tcircle{2}}}
    \pstree[treemode=D, levelsep=0.5cm]{\Tn} {\pstree[levelsep=.8 cm]{\Tc*{3pt}}{\Tcircle{3}\Tcircle{4}}
    }}}
\]
The induced operation $\rho_\sigma:C^\bullet(A,A)^{\otimes 4}\to C^\bullet(A,A^*)$ is given for $f^1, f^2, f^3, f^4\in C^\bullet(A,A)$ and $a_1,\dots,a_{n+1}\in A$ by
\begin{multline*}
\quad\quad\quad \rho_\sigma (f^1,f^2,f^3,f^4)(a_1,\dots,a_n)(a_{n+1}) \\= \sum \pm F(\dots, 1, \dots, D(\dots, f_3(\dots), \dots,f_4(\dots),\dots),\dots)\\
(D(\dots,f_2(\dots),\dots,f_1(\dots),\dots,a_{n+1},\dots)),
\end{multline*}
where the unit $1$ is placed at the special spot of $F$, and the sum is over all possible ways of placing $a_1,\dots, a_{n+1}$ at the dots, such that their cyclic order is preserved and $a_{n+1}$ is at the indicated position.
\item\label{Delta-example}
Our next example will be the $\Delta$-operator. We define $\Delta$ to be the sum of operation associated to the following symbols:
\[\quad
   (-1)\cdot
   \pstree[treemode=L, levelsep=0.8cm, treesep=0.3cm]{\Tr*{1}}
    { \pstree[levelsep=0cm]{\Tc{3pt}}
  { \pstree[treemode=U, levelsep=0.5cm]{\Tn}
     {\Tp[edge={\ncline[linewidth=0.1cm]}]}
    \pstree[treemode=L, levelsep=1cm]{\Tn}
     {\Tcircle{1}}
    \Tn  }}
    +(-1)^{\alpha+\mu+1}
    \pstree[treemode=L, levelsep=0.8cm, treesep=0.3cm]
    {\Tr*{1}} {\pstree{\Tc{3pt}}
    {\pstree[linewidth=0.1cm]{\Tcircle{1}}{\Tp} }}
    +(-1)^{\alpha+\mu+1}
   \pstree[treemode=L, levelsep=0.8cm, treesep=0.3cm]{\Tr*{1}}
    { \pstree[levelsep=0cm]{\Tc{3pt}}
  { \Tn
    \pstree[treemode=L, levelsep=1cm]{\Tn}
     {\Tcircle{1}}
    \pstree[treemode=D, levelsep=0.5cm]{\Tn}
     {\Tp[edge={\ncline[linewidth=0.1cm]}]} }}
\]
Note that its induced operation $C^\bullet(A,A)\to C^\bullet(A,A^*)$ is given for $f\in C^\bullet(A,A)$ and $a_1,\dots, a_{n+1}$ by
\begin{eqnarray*}
\quad\quad\quad \sum \pm F(\dots,a_{n+1},\dots,f(\dots),\dots)(1)
+ \sum \pm F(\dots,f(\dots,a_{n+1},\dots),\dots)(1)\\
+\sum \pm F(\dots,f(\dots),\dots,a_{n+1},\dots)(1),
\end{eqnarray*}
where $f$ is at the special spot of $F$. We see that (up to sign) this sum is $\Delta'$ applied to $F_\sharp(f) \in C^\bullet(A,A^*)$. In order to obtain the correct signs for $\Delta$, we let $\mu=||f||$ denote the degree of $f$, and $\alpha=\sum_{l=1}^{n+1}(|a_l|+1)$ the total degree of the shifted arguments $a_l$. Also, note that $\Delta$ is an operator of degree $-1$, since we apply the unit $1\in (sA)_{+1}$.
\item\label{circ-example}
Consider the symbol
\[
\pstree[treemode=L, levelsep=1cm, treesep=0.3cm, linewidth=0.1cm]
{\Tp} {\pstree[linewidth=0.9pt]{\Tr*{
\begin{pspicture}(0,0)(0.2,0.2) 
 \pscircle[linewidth=0.9pt](0.1,0.1){0.11}
\end{pspicture}}}
{\pstree{\Tcircle{1}} {\Tcircle{2}}}}
\]
Its action on $f^1, f^2\in C^\bullet(A,A)$ and $a_1,\dots,a_{n+1}\in A$ is
$$ \sum \pm F(\dots,f^1(\dots,f^2(\dots),\dots),\dots)(a_{n+1}), $$
where $a_1,\dots,a_n$ are placed at the dots in all possible ways. Note, that (up to sign) we obtain $F_\sharp(f^1\circ f^2)$, where we used the operation ``$\circ$'' from definition \ref{gerst-cup-def}.  With this, we define $[.,.]:C^\bullet(A,A)^{\otimes 2}\to C^\bullet(A,A^*)$ to be the operation associated to
\[
\pstree[treemode=L, levelsep=1cm, treesep=0.3cm, linewidth=0.1cm]
{\Tp} {\pstree[linewidth=0.9pt]{\Tr*{
\begin{pspicture}(0,0)(0.2,0.2) 
 \pscircle[linewidth=0.9pt](0.1,0.1){0.11}
\end{pspicture}}}
{\pstree{\Tcircle{1}} {\Tcircle{2}}}}
-(-1)^{\mu\nu}
\pstree[treemode=L, levelsep=1cm, treesep=0.3cm, linewidth=0.1cm]
{\Tp} {\pstree[linewidth=0.9pt]{\Tr*{
\begin{pspicture}(0,0)(0.2,0.2) 
 \pscircle[linewidth=0.9pt](0.1,0.1){0.11}
\end{pspicture}}}
{\pstree{\Tcircle{2}} {\Tcircle{1}}}}
\]
where $\mu=||f^1||$ and $\nu=||f^2||$ are the degrees of the Hochschild cochains plugged into the bracket.
\item
Define $M:C^\bullet(A,A)^{\otimes 2}\to C^\bullet(A,A^*)$ to be the operation associated to
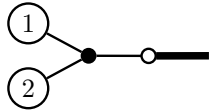
\[
\pstree[treemode=L, levelsep=0.8cm, treesep=0.3cm,
linewidth=0.1cm] {\Tp} {\pstree[linewidth=0.9pt]{\Tr*{
\begin{pspicture}(0,0)(0.2,0.2) 
 \pscircle[linewidth=0.9pt](0.1,0.1){0.11}
\end{pspicture}}}
{\pstree{\Tc*{3pt}} {\Tcircle{1} \Tcircle{2}}}} 
\]
\item\label{unit-reduction-example}
Let $\sigma$ be the following symbol:
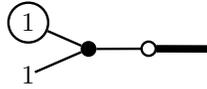
\[
\pstree[treemode=L, levelsep=0.8cm, treesep=0.3cm,
linewidth=0.1cm] {\Tp} {\pstree[linewidth=0.9pt]{\Tr*{
\begin{pspicture}(0,0)(0.2,0.2) 
 \pscircle[linewidth=0.9pt](0.1,0.1){0.11}
\end{pspicture}}}
{\pstree{\Tc*{3pt}} {\Tcircle{1} \Tr*{1}}}} 
\]
Then its action on $f\in C^\bullet(A,A)$ and $a_1,\dots,a_{n+1}\in A$ is given by
\begin{eqnarray*}
\quad & \rho_\sigma(f)(a_1,\dots,a_n)(a_{n+1})& = \sum \pm F(\dots,D(\dots,f(\dots),\dots,1,\dots),\dots)(a_{n+1})\\
&& = \sum\pm F(\dots,f(\dots),\dots)(a_{n+1}),
\end{eqnarray*}
where we have used the fact that $1$ is a strict unit and therefore $D(\dots,f(\dots),\dots,1,\dots)$ is only non-zero, for $D(f(\dots),1)$, which is equal to $\pm f(\dots)$. We thus see that the operation $\rho_\sigma$ is (up to sign) the same as the operation $\rho_{\tilde{\sigma}}$, where  the symbol $\tilde{\sigma}$ is defined by the following graph:
\[
\pstree[treemode=L, levelsep=1cm, treesep=0.3cm, linewidth=0.1cm]
{\Tp} {\pstree[linewidth=0.9pt]{\Tr*{
\begin{pspicture}(0,0)(0.2,0.2) 
 \pscircle[linewidth=0.9pt](0.1,0.1){0.11}
\end{pspicture}}}
{\Tcircle{1}}}
\]
In fact, there are other symbols which also induce the same operations. Here is a third symbol with the same induced operation:
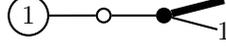
\[
\pstree[treemode=R, levelsep=1.cm, treesep=0.3cm]
  {\Tcircle{1}}{\pstree[levelsep=0.8cm]{\Tc{3pt}}
{\pstree[treemode=R, treesep=0.3cm] {\Tc*{3pt}} {\Tp
[edge={\ncline [linewidth=0.1cm]}] \Tr*{1}}}}
\]
\item
As a final example, consider the following symbol:
 \[
    \pstree[treemode=R, levelsep=1cm, treesep=0.3cm]{\Tcircle{1}}
    { \pstree[levelsep=0cm]{\Tc{3pt}}
  { \pstree[treemode=U, levelsep=0.5cm]{\Tn}
     {\Tp[edge={\ncline[linewidth=0.1cm]}]}
    \pstree[treemode=R, levelsep=1cm]{\Tn}
     {\Tcircle{2}}
    \Tn
  }}
 \]
 Its action on $f^1, f^2\in C^\bullet(A,A)$ and $a_1,\dots, a_{21}\in A$ is 
 $$ \sum\pm F(\dots,a_{21},\dots, f^1(\dots),\dots)(f^2(\dots)). $$
If we assume that the Hochschild cochains $f^k$, for $k=1,2,$ are in components $f^k=\sum_{i\geq 0} f^k_i$, where $f^k_i:sA^{\otimes i}\to sA$, and $F=\sum_{i,j\geq 0} F_{i,j}$, where $F_{i,j}:sA^{\otimes i}\otimes sA\otimes sA^{\otimes j}\to sA^*$, then some of the terms of the above sum with their proper signs are given by
\begin{eqnarray*}
&& (+1)\cdot F_{21,0}(a_1,\dots, a_{21},f^1_0)(f^2_0) \\
&& +(-1)^{(|a_1|+1)\cdot (||f^1||+||f^2||+\sum_{l=2}^{21}(|a_l|+1))}\cdot F_{20,0}(a_2,\dots,a_{21}, f^1_0)(f^2_1(a_1)) \\
&& +\dots \\
&&+(-1)^{\epsilon} F_{6,7}(a_{18},\dots ,a_{21},a_1,a_2,f^1_5(a_3,\dots, a_7),a_8,\dots,a_{14} ) (f^2_3(a_{15},a_{16},a_{17}))\\
&& +\dots,
\end{eqnarray*}
where 
\begin{eqnarray*}
\epsilon&=&\left(\sum_{l=15}^{17}(|a_l|+1)\right)\cdot \left(||f^1||+||f^2||+\sum_{l=18}^{21}(|a_l|+1)\right)\\
&& +\left(\sum_{l=3}^{14}(|a_l|+1)\right)\cdot \left(||f^1||+\sum_{l=18}^{21}(|a_l|+1)\right)\\
&&+\left(\sum_{l=1}^{2}(|a_l|+1)\right)\cdot \left(\sum_{l=18}^{21}(|a_l|+1)\right)
\end{eqnarray*}
\end{enumerate}
\end{expl}
We now define a differential $\delta$ on symbols.
\begin{defn}[Differential $\delta$]\label{differential-definition}
Recall from definition \ref{A-infty-bimodule}, that there are differentials $\delta^A$ and $\delta^{A^*}$ on $C^\bullet(A,A)$ and $C^\bullet(A,A^*)$, respectively. From this, we obtain an induced differential $\delta$ on the space of operations $C^\bullet(A,A)^{\otimes p}\to C^\bullet(A,A^*)$, given for an operation $\rho$, by setting $\delta\rho (f^1,\dots, f^p)$ to be
\begin{equation*}
 \sum_{i=1}^p (-1)^{||\rho||+\sum_{l=1}^{i-1}||f^l||} \rho(f^1,\dots,\delta^A f^i,\dots, f^p) + \delta^{A^*}(\rho(f^1,\dots,f^p)). 
\end{equation*}
Since every symbol $\sigma$ induces an operation $\rho_\sigma:C^\bullet(A,A)^{\otimes p}\to C^\bullet(A,A^*)$, we can define the differential $\delta\sigma$ to be the sum of symbols that represent $\delta \rho_\sigma$. To describe this differential, we need to consider symbols $\sigma'$ together with a chosen vertex $v$ labeled by a filled circle, and a chosen edge $e$ attached to $v$. We denote by $\sigma' / (v,e) $ the symbol obtained by contracting the edge $e$ and replacing $v$ by the other endpoint of $e$. Then, an analysis similar to \cite[section 5]{T} shows, that the differential $\delta\sigma$ consists of a sum $\sum (-1)^{\epsilon_{\sigma'}}\sigma'$ over all $\sigma'$, such that $\sigma=\sigma'/(v,e)$ for some $v$ and $e$. The sign $\epsilon_{\sigma'}$ is obtained by comparing the linear order of the operation associated to $\sigma$ with the one from $\sigma'$. Examples for the differential can be found in the proofs of lemma \ref{2x43} and lemma \ref{XYZ} below.
\end{defn}

We end this subsection with an alternative description for $\Delta$ on cohomology.
\begin{lem}\label{2x43}
$\Delta$ is cohomologous to the operation obtained by the following symbol
\[ 
 (-1)^{\mu}\cdot
    \pstree[treemode=R, levelsep=0.8cm, treesep=0.3cm]{\Tr*{1}}
    { \pstree[levelsep=0cm]{\Tc{3pt}}
  { \pstree[treemode=U, levelsep=0.5cm]{\Tn}
     {\Tcircle{1}}
    \pstree[treemode=R, levelsep=1cm]{\Tn}
     {\Tp[edge={\ncline[linewidth=0.1cm]}]}
    \Tn  }}+
    \pstree[treemode=R, levelsep=0.8cm, treesep=0.3cm]
    {\Tr*{1}} {\pstree{\Tc{3pt}}
    {\pstree[linewidth=0.1cm]{\Tcircle{1}}{\Tp} }}+
    \pstree[treemode=R, levelsep=0.8cm, treesep=0.3cm]{\Tr*{1}}
    { \pstree[levelsep=0cm]{\Tc{3pt}}
  { \Tn
    \pstree[treemode=R, levelsep=1cm]{\Tn}
     {\Tp[edge={\ncline[linewidth=0.1cm]}]}
    \pstree[treemode=D, levelsep=0.5cm]{\Tn}
     {\Tcircle{1}}  }}
\]
where $\mu$ is the degree of the Hochschild cochain plugged into the symbol.
\end{lem}
\begin{proof} 
The proof consists in calculating the differential of a sum of symbols $H$, which will turn out to be the difference of the sum of symbols for $\Delta$ from example \ref{symbol-examples}\eqref{Delta-example} and the sum of symbols stated in the lemma. More precisely, let
\[ H= (-1)\cdot
  \pstree[treemode=L, levelsep=0.8cm, treesep=0.3cm]{\Tr*{1}}
    { \pstree[levelsep=0cm]{\Tc{3pt}}
  { \pstree[treemode=U, levelsep=0.5cm]{\Tn}
     {\Tp[edge={\ncline[linewidth=0.1cm]}] \Tcircle{1}}
    \pstree[treemode=L, levelsep=0.8cm]{\Tn}
     {\Tr*{1}}
    \Tn
  }}-
  \pstree[treemode=L, levelsep=0.8cm, treesep=0.3cm]{\Tr*{1}}
    { \pstree[levelsep=0cm]{\Tc{3pt}}
  { \pstree[treemode=U, levelsep=0.5cm]{\Tn}
     {\pstree[levelsep=0.8cm] {\Tcircle{1}} {\Tn \Tp[edge={\ncline[linewidth=0.1cm]}]}}
    \pstree[treemode=L, levelsep=0.8cm]{\Tn}
     {\Tr*{1}}
    \Tn
  }}-(-1)^{\alpha\cdot\mu}\cdot
  \pstree[treemode=L, levelsep=0.8cm, treesep=0.3cm]{\Tr*{1}}
    { \pstree[levelsep=0cm]{\Tc{3pt}}
  { \pstree[treemode=U, levelsep=0.5cm]{\Tn}
     {\Tcircle{1} \Tp[edge={\ncline[linewidth=0.1cm]}]}
    \pstree[treemode=L, levelsep=0.8cm]{\Tn}
     {\Tr*{1}}
    \Tn
  }},
\]
where $\alpha=\sum_{i=1}^{n+1} (|a_i|+1)$ is the total degree of all input elements $a_1,\dots, a_{n+1}$, which are plugged into the Hochschild cochain. One calculates, that
\begin{multline*} \delta \left( (-1)\cdot
  \pstree[treemode=L, levelsep=0.8cm, treesep=0.3cm]{\Tr*{1}}
    { \pstree[levelsep=0cm]{\Tc{3pt}}
  { \pstree[treemode=U, levelsep=0.5cm]{\Tn}
     {\Tp[edge={\ncline[linewidth=0.1cm]}] \Tcircle{1}}
    \pstree[treemode=L, levelsep=0.8cm]{\Tn}
     {\Tr*{1}}
    \Tn}}
   \right) = (-1)^{\alpha+\mu}\cdot
    \pstree[treemode=L, levelsep=0.8cm, treesep=0.3cm]{\Tr*{1}}
    { \pstree[levelsep=0cm]{\Tc{3pt}}
  { \pstree[treemode=U, levelsep=0.5cm]{\Tn}
     {\Tcircle{1}}
    \pstree[treemode=L, levelsep=1cm]{\Tn}
     {\Tp[edge={\ncline[linewidth=0.1cm]}]}
    \Tn  }} \\+(-1)^{\alpha+\alpha\cdot\mu}\cdot
   \pstree[treemode=R, levelsep=0.8cm, treesep=0.3cm]{\Tr*{1}}
    { \pstree[levelsep=0cm]{\Tc{3pt}}
  { \pstree[treemode=U, levelsep=0.5cm]{\Tn}
     {\Tp[edge={\ncline[linewidth=0.1cm]}]}
    \pstree[treemode=R, levelsep=1cm]{\Tn}
     {\Tcircle{1}}
    \Tn  }}+
 \pstree[treemode=L, levelsep=0.8cm, treesep=0.3cm]{\Tr*{1}}
    { \pstree[levelsep=0cm]{\Tc{3pt}}
  { \pstree[treemode=U, levelsep=0.3cm]{\Tn}
     {\pstree[levelsep=0.5cm]{\Tc*{3pt}}{\Tp[edge={\ncline[linewidth=0.1cm]}]
     \Tcircle{1}}}
    \pstree[treemode=L, levelsep=0.8cm]{\Tn}
     {\Tr*{1}}
    \Tn }},
\end{multline*}
where, in the first two terms, we used the fact described in example \ref{symbol-examples}\eqref{unit-reduction-example}, that the strict unit can only multiply with one term, where it yields the identity. Similarly, we get the other terms as
\begin{multline*} \delta \left( -(-1)^{\alpha\cdot\mu}\cdot
  \pstree[treemode=L, levelsep=0.8cm, treesep=0.3cm]{\Tr*{1}}
    { \pstree[levelsep=0cm]{\Tc{3pt}}
  { \pstree[treemode=U, levelsep=0.5cm]{\Tn}
     {\Tcircle{1} \Tp[edge={\ncline[linewidth=0.1cm]}]}
    \pstree[treemode=L, levelsep=0.8cm]{\Tn}
     {\Tr*{1}}  \Tn}}
   \right) = (-1)^{\alpha\cdot\mu+\alpha+\mu}\cdot
   \pstree[treemode=L, levelsep=0.8cm, treesep=0.3cm]{\Tr*{1}}
    { \pstree[levelsep=0cm]{\Tc{3pt}}
  { \pstree[treemode=U, levelsep=0.5cm]{\Tn}
     {\Tp[edge={\ncline[linewidth=0.1cm]}]}
    \pstree[treemode=L, levelsep=1cm]{\Tn}
     {\Tcircle{1}}
    \Tn  }} \\+(-1)^{\mu}\cdot
    \pstree[treemode=R, levelsep=0.8cm, treesep=0.3cm]{\Tr*{1}}
    { \pstree[levelsep=0cm]{\Tc{3pt}}
  { \pstree[treemode=U, levelsep=0.5cm]{\Tn}
     {\Tcircle{1}}
    \pstree[treemode=R, levelsep=1cm]{\Tn}
     {\Tp[edge={\ncline[linewidth=0.1cm]}]}
    \Tn  }}+(-1)^{\alpha\cdot\mu}\cdot
 \pstree[treemode=L, levelsep=0.8cm, treesep=0.3cm]{\Tr*{1}}
    { \pstree[levelsep=0cm]{\Tc{3pt}}
  { \pstree[treemode=U, levelsep=0.3cm]{\Tn}
     {\pstree[levelsep=0.5cm]{\Tc*{3pt}}
     {\Tcircle{1} \Tp[edge={\ncline[linewidth=0.1cm]}]}}
    \pstree[treemode=L, levelsep=0.8cm]{\Tn}
     {\Tr*{1}}
    \Tn }},
\end{multline*}
\begin{multline*}
   \delta \left( (-1)\cdot
  \pstree[treemode=L, levelsep=0.8cm, treesep=0.3cm]{\Tr*{1}}
    { \pstree[levelsep=0cm]{\Tc{3pt}}
  { \pstree[treemode=U, levelsep=0.5cm]{\Tn}
     {\pstree[levelsep=0.8cm] {\Tcircle{1}} {\Tn \Tp[edge={\ncline[linewidth=0.1cm]}]}}
    \pstree[treemode=L, levelsep=0.8cm]{\Tn}
     {\Tr*{1}}  \Tn }}
   \right) =
    \pstree[treemode=R, levelsep=0.8cm, treesep=0.3cm]
    {\Tr*{1}} {\pstree{\Tc{3pt}}
    {\pstree[linewidth=0.1cm]{\Tcircle{1}}{\Tp} }}
    +(-1)^{\alpha+\mu}\cdot
    \pstree[treemode=L, levelsep=0.8cm, treesep=0.3cm]
    {\Tr*{1}} {\pstree{\Tc{3pt}}
    {\pstree[linewidth=0.1cm]{\Tcircle{1}}{\Tp} }}
    \\
- \pstree[treemode=L, levelsep=0.8cm, treesep=0.3cm]{\Tr*{1}}
    { \pstree[levelsep=0cm]{\Tc{3pt}}
  { \pstree[treemode=U, levelsep=0.3cm]{\Tn}
     {\pstree[levelsep=0.5cm]{\Tc*{3pt}}
     {\Tp[edge={\ncline[linewidth=0.1cm]}] \Tcircle{1}}}
    \pstree[treemode=L, levelsep=0.8cm]{\Tn}
     {\Tr*{1}}
    \Tn }}
    -(-1)^{\alpha\cdot\mu}\cdot
 \pstree[treemode=L, levelsep=0.8cm, treesep=0.3cm]{\Tr*{1}}
    { \pstree[levelsep=0cm]{\Tc{3pt}}
  { \pstree[treemode=U, levelsep=0.3cm]{\Tn}
     {\pstree[levelsep=0.5cm]{\Tc*{3pt}}
     {\Tcircle{1} \Tp[edge={\ncline[linewidth=0.1cm]}]}}
    \pstree[treemode=L, levelsep=0.8cm]{\Tn}
     {\Tr*{1}}
    \Tn }}.
\end{multline*}
By definition \ref{differential-definition} of the differential, we have that $\delta(\rho_H)=\rho_{\delta(H)}$. Thus the operation $\rho_{\delta(H)}$ associated to the sum of symbols in $\delta(H)$ is a boundary in the operation complex of maps $C^\bullet(A,A)\to C^\bullet(A,A^*)$, and must therefore vanishes on Hochschild-cohomology. We see that the following operations are cohomologous, i.e. they induce the same operation on Hochschild-cohomology
\begin{multline*}
\left( (-1)^{\mu}\cdot
    \pstree[treemode=R, levelsep=0.8cm, treesep=0.3cm]{\Tr*{1}}
    { \pstree[levelsep=0cm]{\Tc{3pt}}
  { \pstree[treemode=U, levelsep=0.5cm]{\Tn}
     {\Tcircle{1}}
    \pstree[treemode=R, levelsep=1cm]{\Tn}
     {\Tp[edge={\ncline[linewidth=0.1cm]}]}
    \Tn  }}+
    \pstree[treemode=R, levelsep=0.8cm, treesep=0.3cm]
    {\Tr*{1}} {\pstree{\Tc{3pt}}
    {\pstree[linewidth=0.1cm]{\Tcircle{1}}{\Tp} }}+
    \pstree[treemode=R, levelsep=0.8cm, treesep=0.3cm]{\Tr*{1}}
    { \pstree[levelsep=0cm]{\Tc{3pt}}
  { \Tn
    \pstree[treemode=R, levelsep=1cm]{\Tn}
     {\Tp[edge={\ncline[linewidth=0.1cm]}]}
    \pstree[treemode=D, levelsep=0.5cm]{\Tn}
     {\Tcircle{1}}  }} \right)
  \\ -\left( (-1)^{\alpha\cdot\mu+\alpha+\mu+1}\cdot
   \pstree[treemode=L, levelsep=0.8cm, treesep=0.3cm]{\Tr*{1}}
    { \pstree[levelsep=0cm]{\Tc{3pt}}
  { \pstree[treemode=U, levelsep=0.5cm]{\Tn}
     {\Tp[edge={\ncline[linewidth=0.1cm]}]}
    \pstree[treemode=L, levelsep=1cm]{\Tn}
     {\Tcircle{1}}
    \pstree[treemode=D, levelsep=0.5cm]{\Tn}{\Tn}
      }}
    + (-1)^{\alpha+\mu+1}\cdot
    \pstree[treemode=L, levelsep=0.8cm, treesep=0.3cm]
    {\Tr*{1}} {\pstree{\Tc{3pt}}
    {\pstree[linewidth=0.1cm]{\Tcircle{1}}{\Tp} }}\right. \\
    \left. + (-1)^{\alpha+\mu+1}\cdot
   \pstree[treemode=L, levelsep=0.8cm, treesep=0.3cm]{\Tr*{1}}
    { \pstree[levelsep=0cm]{\Tc{3pt}}
  { \Tn
    \pstree[treemode=L, levelsep=1cm]{\Tn}
     {\Tcircle{1}}
    \pstree[treemode=D, levelsep=0.5cm]{\Tn}
     {\Tp[edge={\ncline[linewidth=0.1cm]}]} }} \right) \cong
       (-1)^{\alpha+\mu}\cdot
    \pstree[treemode=L, levelsep=1cm, treesep=0.3cm]{\Tr*{1}}
    { \pstree[levelsep=0cm]{\Tc{3pt}}
  { \Tn
    \pstree[treemode=L, levelsep=1cm]{\Tn}
     {\Tcircle{1}}
    \pstree[treemode=D, levelsep=0.5cm]{\Tn}
     {\Tp[edge={\ncline[linewidth=0.1cm]}]}
  }}+
    \pstree[treemode=R, levelsep=1cm, treesep=0.3cm]{\Tr*{1}}
    { \pstree[levelsep=0cm]{\Tc{3pt}}
  { \Tn
    \pstree[treemode=R, levelsep=1cm]{\Tn}
     {\Tp[edge={\ncline[linewidth=0.1cm]}]}
    \pstree[treemode=D, levelsep=0.5cm]{\Tn}
     {\Tcircle{1}}
  }}\\+   (-1)^{\alpha\cdot\mu+\mu+1}\cdot
    \pstree[treemode=R, levelsep=1cm, treesep=0.3cm]{\Tr*{1}}
    { \pstree[levelsep=0cm]{\Tc{3pt}}
  { \pstree[treemode=U, levelsep=0.5cm]{\Tn}
     {\Tp[edge={\ncline[linewidth=0.1cm]}]}
    \pstree[treemode=R, levelsep=1cm]{\Tn}
     {\Tcircle{1}}
    \Tn
  }}+   (-1)^{\alpha+\mu+1}\cdot
    \pstree[treemode=L, levelsep=1cm, treesep=0.3cm]{\Tr*{1}}
    { \pstree[levelsep=0cm]{\Tc{3pt}}
  { \pstree[treemode=U, levelsep=0.5cm]{\Tn}
     {\Tcircle{1}}
    \pstree[treemode=L, levelsep=1cm]{\Tn}
     {\Tp[edge={\ncline[linewidth=0.1cm]}]}
    \Tn
  }}.
\end{multline*}
The operations associated to the symbols on the right vanish since our $\infty$-inner product was assumed to be symmetric, and thus is invariant under a $180^\circ$ rotation. The first parenthesis on the left is the sum of symbols stated in the lemma, the second parenthesis is the sum of symbols representing $\Delta$ from example \ref{symbol-examples}\eqref{Delta-example}. This completes the proof of the lemma.
\end{proof}

\subsection{Proof of the BV-relations}\label{proof-bv-section}
This subsection is devoted to the proof of theorem \ref{2x44}. Here, it is crucial, that we can identify $H^\bullet(A,A)$ with $H^\bullet(A,A^*)$ via the map $F_\sharp$ and its quasi-inverse $G_\sharp$ induced by the $\infty$-inner product.

First, notice that $[.,.]$ and $M$ induce a Gerstenhaber-structure on $H^\bullet(A,A)$, since all the homotopies needed in the proof in \cite{G}, and more generally in \cite{GJ2}, still can be applied here. Also, the $\Delta$-operator squares to zero on Hochschild-cohomology, which again can be seen by considering the normalized Hochschild-cochain complex. It is left to show that on Hochschild-cohomology, the deviation of $\Delta$ from being a derivation is given by $[.,.]$. The following symbols $X$, $Y$ and $Z$ will be important ingredients.
\[
X:=
\pstree[treemode=L, levelsep=1cm, treesep=0.3cm, linewidth=0.1cm]
{\Tp} {\pstree[linewidth=0.9pt]{\Tr*{
\begin{pspicture}(0,0)(0.2,0.2) 
 \pscircle[linewidth=0.9pt](0.1,0.1){0.11}
\end{pspicture}}}
{\pstree{\Tcircle{1}} {\Tcircle{2}}}}
\]
\[
Y:= (-1)\cdot
   \pstree[treemode=R, levelsep=0.8cm, treesep=0.3cm,
     linewidth=0.9pt] {\Tr*{1}}{\pstree[levelsep=0cm]{\Tc{3pt}}
  { \pstree[treemode=U, levelsep=0.5cm]{\Tn} {\Tcircle{1}}
    \pstree[treemode=R, levelsep=0.8cm]{\Tn}
       {\pstree{\Tc*{3pt}}{\Tp[edge={\ncline[linewidth=0.1cm]}] \Tcircle{2}}}
    \Tn}}-(-1)\cdot
    \pstree[treemode=R,levelsep=0.8cm, treesep=0.3cm]{\Tr*{1}}
    {\pstree{\Tc{3pt}} {\pstree{\Tcircle{1}}
    {\pstree{\Tc*{3pt}}{\Tp[edge={\ncline[linewidth=0.1cm]}] \Tcircle{2}}}}}
    -(-1)\cdot
\pstree[treemode=R, levelsep=0.8cm, treesep=0.3cm,
     linewidth=0.9pt] {\Tr*{1}}{\pstree[levelsep=0cm]{\Tc{3pt}}
  { \Tn  \pstree[treemode=R, levelsep=0.8cm]{\Tn}
       {\pstree{\Tc*{3pt}} {\Tp[edge={\ncline[linewidth=0.1cm]}] \Tcircle{2}}}
    \pstree[treemode=D, levelsep=0.5cm]{\Tn} {\Tcircle{1}}}}
\]
\[
Z  :=  (-1)\cdot\pstree[treemode=R, levelsep=0.8cm, treesep=0.3cm,
     linewidth=0.9pt] {\Tcircle{2}}{\pstree[levelsep=0cm]{\Tc{3pt}}
  { \pstree[treemode=U, levelsep=0.5cm]{\Tn} {\Tcircle{1}}
    \pstree[treemode=R, levelsep=0.8cm]{\Tn} {\Tp[edge={\ncline[linewidth=0.1cm]}]}
    \Tn}}-(-1)^{\mu\cdot\nu}\cdot
  \pstree[treemode=R, levelsep=0.8cm, treesep=0.3cm,
     linewidth=0.9pt] {\Tcircle{2}}{
    \pstree[treemode=R, levelsep=0.8cm]{\Tc{3pt}}
       {\pstree{\Tcircle{1}}
       {\Tp[edge={\ncline[linewidth=0.1cm]}]}}}
    -(-1)^{\mu\cdot\nu}\cdot
\pstree[treemode=R, levelsep=0.8cm, treesep=0.3cm,
     linewidth=0.9pt] {\Tcircle{2}}{\pstree[levelsep=0cm]{\Tc{3pt}}
  { \Tn
    \pstree[treemode=R, levelsep=0.8cm]{\Tn} {\Tp[edge={\ncline[linewidth=0.1cm]}]}
    \pstree[treemode=D, levelsep=0.5cm]{\Tn} {\Tcircle{1}}}}
\]
\begin{lem} \label{XYZ}
For $f, g\in C^\bullet(A,A)$ with degrees $||f||$ and $||g||$ respectively, the following Hochschild-cochains in $C^\bullet(A,A^*)$ are cohomologous
\begin{enumerate}
\item $X(f,g)\cong Y(f,g)+Z(f,g)$,
\item $Y(f,g)\cong -M(\Delta f, g)$,
\item $Z(f,g)-(-1)^{||f||\cdot||g||}Z(g,f)\cong\Delta(M(f,g))$.
\end{enumerate}
\end{lem}
Before proving lemma \ref{XYZ}, we show how it implies the required BV-relation from theorem \ref{2x44}.
\begin{proof}[Proof of the BV-relation]
From example \ref{symbol-examples}\eqref{circ-example} we see that $[f,g]$ is given as $X(f,g)-(-1)^{||f||\cdot||g||} X(g,f)$. Thus, with lemma \ref{XYZ}, we get
\begin{eqnarray*}
[f,g]&=& X(f,g)-(-1)^{||f||\cdot||g||}X(g,f)\\
&\cong& Y(f,g)-(-1)^{||f||\cdot||g||}Y(g,f)+ Z(f,g)-(-1)^{||f||\cdot||g||}Z(g,f)\\
&\cong& -M(\Delta f,g)+(-1)^{||f||\cdot||g||} M(\Delta g,f)+\Delta(M(f,g))\\
&\cong& \Delta(M(f,g)) - M(\Delta f, g) -(-1)^{||f||} M( f,\Delta g),
\end{eqnarray*}
where we used the graded commutativity of $M$ in the last step.
\end{proof}
It remains to prove lemma \ref{XYZ}.
\begin{proof}[Proof of lemma \ref{XYZ}]
We use the following abbreviation for the relevant degrees. $\alpha=\sum_l (|a_l|+1)$ denotes the total shifted degree of the elements $a_l$ to be applied as arguments, $\mu=||f||$ and $\nu=||g||$, where $f$ and $g$ will be placed into the first, respectively second spot of the corresponding operation of a symbol. 
\begin{enumerate} 
\item\label{H}
Let
\begin{eqnarray*}
H&=& -(-1)^{\mu}\cdot
    \pstree[treemode=R, levelsep=0.8cm, treesep=0.3cm,
     linewidth=0.9pt] {\Tr*{1}}{\pstree[levelsep=0cm]{\Tc{3pt}}
  { \pstree[treemode=U, levelsep=0.5cm]{\Tn} {\Tcircle{1}}
    \pstree[treemode=R, levelsep=0.8cm]{\Tn} {\Tp[edge={\ncline[linewidth=0.1cm]}]}
    \pstree[treemode=D, levelsep=0.5cm]{\Tn} {\Tcircle{2}}}}
    -(-1)^{\mu\cdot\nu}\cdot
    \pstree[treemode=R, levelsep=0.8cm, treesep=0.3cm,
     linewidth=0.9pt] {\Tr*{1}}{\pstree[levelsep=0cm]{\Tc{3pt}}
  { \Tn
    \pstree[treemode=R, levelsep=0.8cm]{\Tn}
        {\pstree{\Tcircle{1}}{\Tp[edge={\ncline[linewidth=0.1cm]}]}}
    \pstree[treemode=D, levelsep=0.5cm]{\Tn} {\Tcircle{2}}}}
    -
    \pstree[treemode=R,levelsep=0.8cm, treesep=0.3cm]{\Tr*{1}}
    {\pstree{\Tc{3pt}} {\pstree{\Tcircle{1}}
    {\Tp[edge={\ncline[linewidth=0.1cm]}] \Tcircle{2}}}} \\
&& -(-1)^{\mu\cdot\nu}\cdot
    \pstree[treemode=R, levelsep=0.8cm, treesep=0.3cm,
     linewidth=0.9pt] {\Tr*{1}}{\pstree[levelsep=0cm]{\Tc{3pt}}
  { \Tn
    \pstree[treemode=R, levelsep=0.8cm]{\Tn} {\Tp[edge={\ncline[linewidth=0.1cm]}]}
    \pstree[treemode=D, levelsep=0.5cm]{\Tn} {\Tcircle{2} \Tcircle{1}}}}-
   \pstree[treemode=R, levelsep=0.8cm, treesep=0.3cm,
     linewidth=0.9pt] {\Tr*{1}}{\pstree[levelsep=0cm]{\Tc{3pt}}
  { \Tn
    \pstree[treemode=R, levelsep=0.8cm]{\Tn} {\Tp[edge={\ncline[linewidth=0.1cm]}]}
    \pstree[treemode=D, levelsep=0.5cm]{\Tn}
       {\pstree[levelsep=1cm]{\Tcircle{1}} {\Tcircle{2}}}}}-
   \pstree[treemode=R, levelsep=0.8cm, treesep=0.3cm,
     linewidth=0.9pt] {\Tr*{1}}{\pstree[levelsep=0cm]{\Tc{3pt}}
  { \Tn
    \pstree[treemode=R, levelsep=0.8cm]{\Tn} {\Tp[edge={\ncline[linewidth=0.1cm]}]}
    \pstree[treemode=D, levelsep=0.5cm]{\Tn} {\Tcircle{1}
    \Tcircle{2}}}}.
\end{eqnarray*}
We will show that the differential of $H$ is exactly $\delta(H)=-X+Y+Z$. Before doing so, we want to remark on the motivation for $H$. As it was mentioned in the introduction, there is a close connection to the Chas-Sullivan BV-structure from \cite{CS}. In this analogy, $H$ should be compared to the homotopy of \cite[figure 7]{CS}.

We now calculate the boundary $\delta(H)$ term by term:
\begin{multline*} 
    \delta \left( -(-1)^{\mu}\cdot
     \pstree[treemode=R, levelsep=0.8cm, treesep=0.3cm,
     linewidth=0.9pt] {\Tr*{1}}{\pstree[levelsep=0cm]{\Tc{3pt}}
  { \pstree[treemode=U, levelsep=0.5cm]{\Tn} {\Tcircle{1}}
    \pstree[treemode=R, levelsep=0.8cm]{\Tn} {\Tp[edge={\ncline[linewidth=0.1cm]}]}
    \pstree[treemode=D, levelsep=0.5cm]{\Tn} {\Tcircle{2}}}}\right)  =
(-1)\cdot \pstree[treemode=R, levelsep=0.8cm, treesep=0.3cm,
     linewidth=0.9pt] {\Tcircle{2}}{\pstree[levelsep=0cm]{\Tc{3pt}}
  { \pstree[treemode=U, levelsep=0.5cm]{\Tn} {\Tcircle{1}}
    \pstree[treemode=R, levelsep=0.8cm]{\Tn} {\Tp[edge={\ncline[linewidth=0.1cm]}]}
    \Tn}}
    -(-1)^{\alpha\cdot\mu+\mu\cdot\nu}\cdot
    \pstree[treemode=R,levelsep=0.8cm, treesep=0.3cm]{\Tr*{1}}
    {\pstree{\Tc{3pt}} {\pstree{\Tc*{3pt}}
    {\Tcircle{1} \Tp[edge={\ncline[linewidth=0.1cm]}]
    \Tcircle{2}}}} \\
+ \pstree[treemode=R, levelsep=0.8cm, treesep=0.3cm,
     linewidth=0.9pt] {\Tcircle{1}}{\pstree[levelsep=0cm]{\Tc{3pt}}
  { \Tn
    \pstree[treemode=R, levelsep=0.8cm]{\Tn} {\Tp[edge={\ncline[linewidth=0.1cm]}]}
    \pstree[treemode=D, levelsep=0.5cm]{\Tn} {\Tcircle{2}}}}
 -(-1)^{\alpha\cdot\mu+\mu\cdot\nu+\nu}\cdot
    \pstree[treemode=R, levelsep=0.8cm, treesep=0.3cm,
     linewidth=0.9pt] {\Tr*{1}}{\pstree[levelsep=0cm]{\Tc{3pt}}
  { \Tn  \pstree[treemode=R, levelsep=0.8cm]{\Tn}
       {\pstree{\Tc*{3pt}} {\Tcircle{1} \Tp[edge={\ncline[linewidth=0.1cm]}]}}
    \pstree[treemode=D, levelsep=0.5cm]{\Tn} {\Tcircle{2}}}}-
   \pstree[treemode=R, levelsep=0.8cm, treesep=0.3cm,
     linewidth=0.9pt] {\Tr*{1}}{\pstree[levelsep=0cm]{\Tc{3pt}}
  { \pstree[treemode=U, levelsep=0.5cm]{\Tn} {\Tcircle{1}}
    \pstree[treemode=R, levelsep=0.8cm]{\Tn}
       {\pstree{\Tc*{3pt}}{\Tp[edge={\ncline[linewidth=0.1cm]}] \Tcircle{2}}}
    \Tn}}
\end{multline*}
\begin{multline*} 
    \delta \left( -(-1)^{\mu\cdot\nu}\cdot
   \pstree[treemode=R, levelsep=0.8cm, treesep=0.3cm,
     linewidth=0.9pt] {\Tr*{1}}{\pstree[levelsep=0cm]{\Tc{3pt}}
  { \Tn  \pstree[treemode=R, levelsep=0.8cm]{\Tn}
          {\pstree{\Tcircle{1}}{\Tp[edge={\ncline[linewidth=0.1cm]}]}}
    \pstree[treemode=D, levelsep=0.5cm]{\Tn} {\Tcircle{2}}}}\right)  =
   (-1)^{\alpha\cdot\mu+\mu\cdot\nu+\nu}\cdot
    \pstree[treemode=R, levelsep=0.8cm, treesep=0.3cm,
     linewidth=0.9pt] {\Tr*{1}}{\pstree[levelsep=0cm]{\Tc{3pt}}
  { \Tn  \pstree[treemode=R, levelsep=0.8cm]{\Tn}
       {\pstree{\Tc*{3pt}} {\Tcircle{1} \Tp[edge={\ncline[linewidth=0.1cm]}]}}
    \pstree[treemode=D, levelsep=0.5cm]{\Tn} {\Tcircle{2}}}} \\
   +(-1)^{\mu\cdot\nu+\nu}\cdot
   \pstree[treemode=R, levelsep=0.8cm, treesep=0.3cm,
     linewidth=0.9pt] {\Tr*{1}}{\pstree[levelsep=0cm]{\Tc{3pt}}
  { \Tn  \pstree[treemode=R, levelsep=0.8cm]{\Tn}
       {\pstree{\Tc*{3pt}} {\Tp[edge={\ncline[linewidth=0.1cm]}] \Tcircle{1}}}
    \pstree[treemode=D, levelsep=0.5cm]{\Tn} {\Tcircle{2}}}}
  -(-1)^{\mu\cdot\nu}\cdot
   \pstree[treemode=R, levelsep=0.8cm, treesep=0.3cm,
     linewidth=0.9pt] {\Tcircle{2}}{
    \pstree[treemode=R, levelsep=0.8cm]{\Tc{3pt}}
       {\pstree{\Tcircle{1}}
       {\Tp[edge={\ncline[linewidth=0.1cm]}]}}} \\
 -(-1)^{\mu\cdot\nu}\cdot
 \pstree[treemode=R, levelsep=0.8cm, treesep=0.3cm,
     linewidth=0.9pt] {\Tr*{1}}{\pstree[treemode=R, levelsep=0.8cm]{\Tc{3pt}}
       {\pstree{\Tc*{3pt}} {\pstree{\Tcircle{1}}{\Tp[edge={\ncline[linewidth=0.1cm]}]}
         \Tcircle{2}}}}
\end{multline*}
\begin{multline*} 
    \delta \left( (-1)\cdot
    \pstree[treemode=R,levelsep=0.8cm, treesep=0.3cm]{\Tr*{1}}
    {\pstree{\Tc{3pt}} {\pstree{\Tcircle{1}}
    {\Tp[edge={\ncline[linewidth=0.1cm]}] \Tcircle{2}}}}
    \right)  =
    (-1)^{\alpha\cdot\mu+\mu\cdot\nu}\cdot
    \pstree[treemode=R,levelsep=0.8cm, treesep=0.3cm]{\Tr*{1}}
    {\pstree{\Tc{3pt}} {\pstree{\Tc*{3pt}}
    {\Tcircle{1} \Tp[edge={\ncline[linewidth=0.1cm]}]
    \Tcircle{2}}}} \\
    +(-1)^{\mu\cdot\nu}\cdot
    \pstree[treemode=R,levelsep=0.8cm, treesep=0.3cm]{\Tr*{1}}
    {\pstree{\Tc{3pt}} {\pstree{\Tc*{3pt}}
    { \Tp[edge={\ncline[linewidth=0.1cm]}] \Tcircle{1}
    \Tcircle{2}}}}
   +\pstree[treemode=R,levelsep=0.8cm, treesep=0.3cm]{\Tr*{1}}
    {\pstree{\Tc{3pt}} {\pstree{\Tc*{3pt}}
    {\Tp[edge={\ncline[linewidth=0.1cm]}]
    \Tcircle{2} \Tcircle{1}}}}
   -(-1)^{\mu}\cdot
    \pstree[treemode=R,levelsep=0.8cm, treesep=0.3cm]{\Tr*{1}}
    {\pstree{\Tc{3pt}} {\pstree{\Tcircle{1}}
    {\pstree{\Tc*{3pt}}{\Tp[edge={\ncline[linewidth=0.1cm]}] \Tcircle{2}}}}} \\
  +(-1)^{\mu\cdot\nu}\cdot
   \pstree[treemode=R,levelsep=0.8cm, treesep=0.3cm]{\Tr*{1}}
    {\pstree{\Tc{3pt}} {\pstree{\Tc*{3pt}}
    {\pstree{\Tcircle{1}}{\Tp[edge={\ncline[linewidth=0.1cm]}]} \Tcircle{2}}}}+
    \pstree[treemode=R,levelsep=0.8cm, treesep=0.3cm]{\Tr*{1}}
    {\pstree{\Tc{3pt}} {\pstree{\Tc*{3pt}}
    {\Tp[edge={\ncline[linewidth=0.1cm]}]
    \pstree{\Tcircle{1}}{\Tcircle{2}}}}}
\end{multline*}
\begin{multline*} 
    \delta \left( -(-1)^{\mu\cdot\nu}\cdot
   \pstree[treemode=R, levelsep=0.8cm, treesep=0.3cm,
     linewidth=0.9pt] {\Tr*{1}}{\pstree[levelsep=0cm]{\Tc{3pt}}
  { \Tn
    \pstree[treemode=R, levelsep=0.8cm]{\Tn} {\Tp[edge={\ncline[linewidth=0.1cm]}]}
    \pstree[treemode=D, levelsep=0.5cm]{\Tn} {\Tcircle{2} \Tcircle{1}}}}
    \right)  =
   -(-1)^{\mu\cdot\nu}\cdot
    \pstree[treemode=R, levelsep=0.8cm, treesep=0.3cm,
     linewidth=0.9pt] {\Tcircle{2}}{\pstree[levelsep=0cm]{\Tc{3pt}}
  { \Tn
    \pstree[treemode=R, levelsep=0.8cm]{\Tn} {\Tp[edge={\ncline[linewidth=0.1cm]}]}
    \pstree[treemode=D, levelsep=0.5cm]{\Tn} {\Tcircle{1}}}}
  -(-1)^{\mu\cdot\nu}\cdot
   \pstree[treemode=R, levelsep=0.8cm, treesep=0.3cm,
     linewidth=0.9pt] {\Tr*{1}}{\pstree[levelsep=0cm]{\Tc{3pt}}
   { \Tn
    \pstree[treemode=R, levelsep=0.8cm]{\Tn} {\Tp[edge={\ncline[linewidth=0.1cm]}]}
    \pstree[treemode=D, levelsep=0.4cm]{\Tn} {
    \pstree[levelsep=0.7cm]{\Tc*{3pt}}{\Tcircle{2}
    \Tcircle{1}}}}} \\
  -(-1)^{\mu\cdot\nu+\nu}\cdot
  \pstree[treemode=R, levelsep=0.8cm, treesep=0.3cm,
     linewidth=0.9pt] {\Tr*{1}}{\pstree[levelsep=0cm]{\Tc{3pt}}
  { \Tn  \pstree[treemode=R, levelsep=0.8cm]{\Tn}
       {\pstree{\Tc*{3pt}} {\Tp[edge={\ncline[linewidth=0.1cm]}] \Tcircle{1}}}
    \pstree[treemode=D, levelsep=0.5cm]{\Tn} {\Tcircle{2}}}}
   -(-1)^{\mu\cdot\nu}\cdot
    \pstree[treemode=R,levelsep=0.8cm, treesep=0.3cm]{\Tr*{1}}
    {\pstree{\Tc{3pt}} {\pstree{\Tc*{3pt}}
    { \Tp[edge={\ncline[linewidth=0.1cm]}] \Tcircle{1}
    \Tcircle{2}}}}
\end{multline*}
\begin{multline*} 
    \delta \left( (-1)\cdot
   \pstree[treemode=R, levelsep=0.8cm, treesep=0.3cm,
     linewidth=0.9pt] {\Tr*{1}}{\pstree[levelsep=0cm]{\Tc{3pt}}
  { \Tn
    \pstree[treemode=R, levelsep=0.8cm]{\Tn} {\Tp[edge={\ncline[linewidth=0.1cm]}]}
    \pstree[treemode=D, levelsep=0.5cm]{\Tn}
       {\pstree[levelsep=1cm]{\Tcircle{1}} {\Tcircle{2}}}}}
    \right)  =
    (-1)\cdot
    \pstree[treemode=R, levelsep=1cm, treesep=0.3cm,
     linewidth=0.9pt] {\Tcircle{2}}{
    \pstree[treemode=R, levelsep=0.8cm]{\Tcircle{1}}
       {\pstree{\Tc{3pt}}
       {\Tp[edge={\ncline[linewidth=0.1cm]}]}}}-
    \pstree[treemode=R,levelsep=0.8cm, treesep=0.3cm]{\Tr*{1}}
    {\pstree{\Tc{3pt}} {\pstree{\Tc*{3pt}}
    {\Tp[edge={\ncline[linewidth=0.1cm]}]
    \pstree{\Tcircle{1}}{\Tcircle{2}}}}} \\
 + \pstree[treemode=R, levelsep=0.8cm, treesep=0.3cm,
     linewidth=0.9pt] {\Tr*{1}}{\pstree[levelsep=0cm]{\Tc{3pt}}
  { \Tn
    \pstree[treemode=R, levelsep=0.8cm]{\Tn} {\Tp[edge={\ncline[linewidth=0.1cm]}]}
    \pstree[treemode=D, levelsep=0.4cm]{\Tn} {
    \pstree[levelsep=0.7cm]{\Tc*{3pt}}{\Tcircle{1}
    \Tcircle{2}}}}} + (-1)^{\mu\cdot\nu}\cdot
  \pstree[treemode=R, levelsep=0.8cm, treesep=0.3cm,
     linewidth=0.9pt] {\Tr*{1}}{\pstree[levelsep=0cm]{\Tc{3pt}}
  { \Tn
    \pstree[treemode=R, levelsep=0.8cm]{\Tn} {\Tp[edge={\ncline[linewidth=0.1cm]}]}
    \pstree[treemode=D, levelsep=0.4cm]{\Tn} {
    \pstree[levelsep=0.7cm]{\Tc*{3pt}}{\Tcircle{2}
    \Tcircle{1}}}}}
\end{multline*}
\begin{multline*} 
    \delta \left( (-1)\cdot
   \pstree[treemode=R, levelsep=0.8cm, treesep=0.3cm,
     linewidth=0.9pt] {\Tr*{1}}{\pstree[levelsep=0cm]{\Tc{3pt}}
  { \Tn
    \pstree[treemode=R, levelsep=0.8cm]{\Tn} {\Tp[edge={\ncline[linewidth=0.1cm]}]}
    \pstree[treemode=D, levelsep=0.5cm]{\Tn} {\Tcircle{1} \Tcircle{2}}}}
    \right)  =
  (-1)\cdot\pstree[treemode=R, levelsep=0.8cm, treesep=0.3cm,
     linewidth=0.9pt] {\Tcircle{1}}{\pstree[levelsep=0cm]{\Tc{3pt}}
  { \Tn
    \pstree[treemode=R, levelsep=0.8cm]{\Tn} {\Tp[edge={\ncline[linewidth=0.1cm]}]}
    \pstree[treemode=D, levelsep=0.5cm]{\Tn} {\Tcircle{2}}}}-
  \pstree[treemode=R, levelsep=0.8cm, treesep=0.3cm,
     linewidth=0.9pt] {\Tr*{1}}{\pstree[levelsep=0cm]{\Tc{3pt}}
  { \Tn
    \pstree[treemode=R, levelsep=0.8cm]{\Tn} {\Tp[edge={\ncline[linewidth=0.1cm]}]}
    \pstree[treemode=D, levelsep=0.4cm]{\Tn} {
    \pstree[levelsep=0.7cm]{\Tc*{3pt}}{\Tcircle{1}
    \Tcircle{2}}}}} \\
  -(-1)\cdot \pstree[treemode=R, levelsep=0.8cm, treesep=0.3cm,
     linewidth=0.9pt] {\Tr*{1}}{\pstree[levelsep=0cm]{\Tc{3pt}}
  { \Tn  \pstree[treemode=R, levelsep=0.8cm]{\Tn}
       {\pstree{\Tc*{3pt}} {\Tp[edge={\ncline[linewidth=0.1cm]}] \Tcircle{2}}}
    \pstree[treemode=D, levelsep=0.5cm]{\Tn} {\Tcircle{1}}}}-
    \pstree[treemode=R,levelsep=0.8cm, treesep=0.3cm]{\Tr*{1}}
    {\pstree{\Tc{3pt}} {\pstree{\Tc*{3pt}}
    { \Tp[edge={\ncline[linewidth=0.1cm]}] \Tcircle{2}
    \Tcircle{1}}}}
\end{multline*}
A close investigation shows that all terms in $\delta(H)$ cancel, except for the terms in $-X+Y+Z$. This implies the claim.
\item
In order to see that $Y(f,g)\cong -M(\Delta f, g)$, we consider the following differential:
\begin{equation*}
\delta \left(\,\, (-1)\cdot\pstree[treemode=R,
levelsep=1.2cm, treesep=0.3cm, linewidth=0.9pt]
{\Tr*{\begin{pspicture}(0,0.15)(0.5,0.5)
 \rput[b](0.25,0.17){$1$} \pscircle(0.25,0.25){0.3}
\end{pspicture}}}{\pstree[levelsep=0cm]{\Tc{3pt}}
  { \Tn \pstree[treemode=R, levelsep=1cm]{\Tn}
    {\Tp[edge={\ncline[linewidth=0.1cm]}]}
    \pstree[treemode=D, levelsep=0.7cm]{\Tn} {\Tcircle{2}}}} \,\,\right) =
\pstree[treemode=L, levelsep=0.8cm, treesep=0.3cm,
linewidth=0.1cm] {\Tp} {\pstree[linewidth=0.9pt]{\Tr*{
\begin{pspicture}(0,0)(0.2,0.2) 
 \pscircle[linewidth=0.9pt](0.1,0.1){0.11}
\end{pspicture}}}
{\pstree{\Tc*{3pt}} {\Tr*{
\begin{pspicture}(0,0.15)(0.5,0.5)
 \rput[b](0.25,0.17){$1$}
 \pscircle(0.25,0.25){0.3}
\end{pspicture}
} \Tcircle{2}}}} -(-1)^{\mu}\cdot
\pstree[treemode=R, levelsep=1.2cm, treesep=0.3cm]
  {\Tr*{
\begin{pspicture}(0,0.15)(0.5,0.5)
 \rput[b](0.25,0.17){$1$}
 \pscircle(0.25,0.25){0.3}
\end{pspicture} }}{\pstree[levelsep=0.8cm]{\Tc{3pt}}
{\pstree[treemode=R, treesep=0.3cm] {\Tc*{3pt}} {\Tp
[edge={\ncline [linewidth=0.1cm]}] \Tcircle{2}}}}
\end{equation*}
Since $\Delta$ is a chain map, we can pre-compose this equation by $\Delta\otimes \id$, which shows that the following induced operations are cohomologous,
\[
\pstree[treemode=L, levelsep=0.8cm, treesep=0.3cm,
linewidth=0.1cm] {\Tp} {\pstree[linewidth=0.9pt]{\Tr*{
\begin{pspicture}(0,0)(0.2,0.2) 
 \pscircle[linewidth=0.9pt](0.1,0.1){0.11}
\end{pspicture}}}
{\pstree{\Tc*{3pt}} {\Tr*{\Delta(\,
\begin{pspicture}(0,0.15)(0.5,0.5)
 \rput[b](0.25,0.17){$1$}
 \pscircle(0.25,0.25){0.3}
\end{pspicture}
\,)} \Tcircle{2}}}} \cong(-1)^{\mu}\cdot
\pstree[treemode=R, levelsep=1.2cm, treesep=0.3cm]
  {\Tr*{\Delta(\,
\begin{pspicture}(0,0.15)(0.5,0.5)
 \rput[b](0.25,0.17){$1$}
 \pscircle(0.25,0.25){0.3}
\end{pspicture} \,)}}{\pstree[levelsep=0.8cm]{\Tc{3pt}}
{\pstree[treemode=R, treesep=0.3cm] {\Tc*{3pt}} {\Tp
[edge={\ncline [linewidth=0.1cm]}] \Tcircle{2}}}}.
\]
The left term is exactly $F\circ( M(\Delta f, g))$, when applying $f$ and $g$ to it. The term on the right will be seen to be equal to $Y(f,g)$, after applying the description for $\Delta$ from lemma \ref{2x43}:
\begin{multline*}
\quad\quad\quad (-1)^{\mu}\cdot \pstree[treemode=R, levelsep=1.2cm,
treesep=0.3cm]
  {\Tr*{\Delta(\,
\begin{pspicture}(0,0.15)(0.5,0.5)
 \rput[b](0.25,0.17){$1$}
 \pscircle(0.25,0.25){0.3}
\end{pspicture} \,)}}{\pstree[levelsep=0.8cm]{\Tc{3pt}}
{\pstree[treemode=R, treesep=0.3cm] {\Tc*{3pt}} {\Tp
[edge={\ncline [linewidth=0.1cm]}] \Tcircle{2}}}} \cong (-1)\cdot
   \pstree[treemode=R, levelsep=0.8cm, treesep=0.3cm,
     linewidth=0.9pt] {\Tr*{1}}{\pstree[levelsep=0cm]{\Tc{3pt}}
  { \pstree[treemode=U, levelsep=0.5cm]{\Tn} {\Tcircle{1}}
    \pstree[treemode=R, levelsep=0.8cm]{\Tn}
       {\pstree{\Tc*{3pt}}{\Tp[edge={\ncline[linewidth=0.1cm]}] \Tcircle{2}}}
    \Tn}} \\-(-1)^{\mu}\cdot
    \pstree[treemode=R,levelsep=0.8cm, treesep=0.3cm]{\Tr*{1}}
    {\pstree{\Tc{3pt}} {\pstree{\Tcircle{1}}
    {\pstree{\Tc*{3pt}}{\Tp[edge={\ncline[linewidth=0.1cm]}] \Tcircle{2}}}}}
-(-1)^{\mu}\cdot \pstree[treemode=R, levelsep=0.8cm,
treesep=0.3cm,
     linewidth=0.9pt] {\Tr*{1}}{\pstree[levelsep=0cm]{\Tc{3pt}}
  { \Tn  \pstree[treemode=R, levelsep=0.8cm]{\Tn}
       {\pstree{\Tc*{3pt}} {\Tp[edge={\ncline[linewidth=0.1cm]}] \Tcircle{2}}}
    \pstree[treemode=D, levelsep=0.5cm]{\Tn} {\Tcircle{1}}}}
\end{multline*}
\item
We now show that $\Delta(M(f,g))\cong Z(f,g)-(-1)^{||f|| \cdot ||g||}\cdot Z(g,f)$. For this, we calculate the differential of 
\begin{multline*}
H'= -(-1)^{\mu}\cdot
    \pstree[treemode=R, levelsep=0.8cm,treesep=0.3cm,
     linewidth=0.9pt] {\Tr*{1}}{\pstree[levelsep=0cm]{\Tc{3pt}}
  { \pstree[treemode=U, levelsep=0.5cm]{\Tn} {\Tcircle{1}}
    \pstree[treemode=R, levelsep=0.8cm]{\Tn} {\Tp[edge={\ncline[linewidth=0.1cm]}]}
    \pstree[treemode=D, levelsep=0.5cm]{\Tn} {\Tcircle{2}}}}
    -(-1)^{\mu\cdot\nu}\cdot
    \pstree[treemode=R, levelsep=0.8cm, treesep=0.3cm,
     linewidth=0.9pt] {\Tr*{1}}{\pstree[levelsep=0cm]{\Tc{3pt}}
  { \Tn
    \pstree[treemode=R, levelsep=0.8cm]{\Tn}
          {\pstree{\Tcircle{1}}{\Tp[edge={\ncline[linewidth=0.1cm]}]}}
    \pstree[treemode=D, levelsep=0.5cm]{\Tn} {\Tcircle{2}}}} \\
  -(-1)^{\mu\cdot\nu}\cdot
    \pstree[treemode=R, levelsep=0.8cm, treesep=0.3cm,
     linewidth=0.9pt] {\Tr*{1}}{\pstree[levelsep=0cm]{\Tc{3pt}}
  { \Tn
    \pstree[treemode=R, levelsep=0.8cm]{\Tn} {\Tp[edge={\ncline[linewidth=0.1cm]}]}
    \pstree[treemode=D, levelsep=0.5cm]{\Tn} {\Tcircle{2} \Tcircle{1}}}}
  -(-1)^{\mu}\cdot
   \pstree[treemode=R,levelsep=0.8cm, treesep=0.3cm,
     linewidth=0.9pt] {\Tr*{1}}{\pstree[levelsep=0cm]{\Tc{3pt}}
  { \pstree[treemode=U, levelsep=0.5cm]{\Tn} {\Tcircle{1}}
    \pstree[treemode=R, levelsep=0.8cm]{\Tn}
          {\pstree{\Tcircle{2}}{\Tp[edge={\ncline[linewidth=0.1cm]}]}}
    \Tn}}
  -(-1)^{\mu+\nu+\mu\cdot\nu}\cdot
    \pstree[treemode=R, levelsep=0.8cm, treesep=0.3cm,
     linewidth=0.9pt] {\Tr*{1}}{\pstree[levelsep=0cm]{\Tc{3pt}}
  { \pstree[treemode=U, levelsep=0.5cm]{\Tn} {\Tcircle{1} \Tcircle{2}}
    \pstree[treemode=R, levelsep=0.8cm]{\Tn} {\Tp[edge={\ncline[linewidth=0.1cm]}]}
    \Tn}}.
\end{multline*}
We remark, that the first three terms were already used in part \eqref{H} as terms in $H$, and will now be subtracted again. We chose this way of proceeding, in order to have a closer analogy to the proof in \cite{CS}. So, for the differential of $H'$, we will only perform the calculation for the last two terms, and refer to the calculation in \eqref{H} for the first three terms.
\begin{multline*} 
    \delta\left( -(-1)^{\mu}\cdot
     \pstree[treemode=R,levelsep=0.8cm, treesep=0.3cm,
     linewidth=0.9pt] {\Tr*{1}}{\pstree[levelsep=0cm]{\Tc{3pt}}
  { \pstree[treemode=U, levelsep=0.5cm]{\Tn} {\Tcircle{1}}
    \pstree[treemode=R, levelsep=0.8cm]{\Tn}
          {\pstree{\Tcircle{2}}{\Tp[edge={\ncline[linewidth=0.1cm]}]}}
    \Tn}}\right)  \cong
    (-1)^{\alpha\cdot\nu}\cdot
     \pstree[treemode=R,levelsep=0.8cm,treesep=0.3cm,
     linewidth=0.9pt] {\Tr*{1}}{\pstree[levelsep=0cm]{\Tc{3pt}}
  { \pstree[treemode=U, levelsep=0.5cm]{\Tn} {\Tcircle{1}}
    \pstree[treemode=R, levelsep=0.8cm]{\Tn}
  {\pstree{\Tc*{3pt}}{\Tcircle{2}\Tp[edge={\ncline[linewidth=0.1cm]}]}} \Tn}}+
\pstree[treemode=R,levelsep=0.8cm, treesep=0.3cm,
     linewidth=0.9pt] {\Tr*{1}}{\pstree[levelsep=0cm]{\Tc{3pt}}
  { \pstree[treemode=U, levelsep=0.5cm]{\Tn} {\Tcircle{1}}
    \pstree[treemode=R, levelsep=0.8cm]{\Tn}
  {\pstree{\Tc*{3pt}}{\Tp[edge={\ncline[linewidth=0.1cm]}]\Tcircle{2}}} \Tn}}\\
 +\pstree[treemode=R, levelsep=0.8cm, treesep=0.3cm,
     linewidth=0.9pt] {\Tcircle{1}}{
    \pstree[treemode=R, levelsep=0.8cm]{\Tc{3pt}}
       {\pstree{\Tcircle{2}}
       {\Tp[edge={\ncline[linewidth=0.1cm]}]}}}
  -(-1)^{\alpha\cdot\mu+\mu\cdot\nu}\cdot
  \pstree[treemode=R, levelsep=0.8cm, treesep=0.3cm,
     linewidth=0.9pt] {\Tr*{1}}{\pstree[treemode=R, levelsep=0.8cm]{\Tc{3pt}}
       {\pstree{\Tc*{3pt}} {\Tcircle{1}
       \pstree{\Tcircle{2}}{\Tp[edge={\ncline[linewidth=0.1cm]}]} }}}
\end{multline*}
\begin{multline*} 
    \delta\left( -(-1)^{\mu+\nu+\mu\cdot\nu}\cdot
    \pstree[treemode=R, levelsep=0.8cm, treesep=0.3cm,
     linewidth=0.9pt] {\Tr*{1}}{\pstree[levelsep=0cm]{\Tc{3pt}}
  { \pstree[treemode=U, levelsep=0.5cm]{\Tn} {\Tcircle{1} \Tcircle{2}}
    \pstree[treemode=R, levelsep=0.8cm]{\Tn} {\Tp[edge={\ncline[linewidth=0.1cm]}]}
    \Tn}}    \right)  \cong
   (-1)^{\mu\cdot\nu}\cdot
   \pstree[treemode=R, levelsep=0.8cm, treesep=0.3cm,
     linewidth=0.9pt] {\Tcircle{1}}{\pstree[levelsep=0cm]{\Tc{3pt}}
  { \pstree[treemode=U, levelsep=0.5cm]{\Tn} {\Tcircle{2}}
    \pstree[treemode=R, levelsep=0.8cm]{\Tn} {\Tp[edge={\ncline[linewidth=0.1cm]}]}
    \Tn}}+
  (-1)^{\mu+\nu+\mu\cdot\nu}\cdot
    \pstree[treemode=R, levelsep=0.8cm, treesep=0.3cm,
     linewidth=0.9pt] {\Tr*{1}}{\pstree[levelsep=0cm]{\Tc{3pt}}
  { \pstree[treemode=U, levelsep=0.4cm]{\Tn} {
    \pstree[levelsep=0.7cm]{\Tc*{3pt}}{\Tcircle{1}
    \Tcircle{2}}}
    \pstree[treemode=R, levelsep=0.8cm]{\Tn}
    {\Tp[edge={\ncline[linewidth=0.1cm]}]}\Tn}} \\
  -(-1)^{\alpha\cdot\nu}\cdot
   \pstree[treemode=R, levelsep=0.8cm, treesep=0.3cm,
     linewidth=0.9pt] {\Tr*{1}}{\pstree[levelsep=0cm]{\Tc{3pt}}
  { \pstree[treemode=U, levelsep=0.5cm]{\Tn} {\Tcircle{1}}
    \pstree[treemode=R, levelsep=0.8cm]{\Tn}
       {\pstree{\Tc*{3pt}} {\Tcircle{2}
                            \Tp[edge={\ncline[linewidth=0.1cm]}] }}
    \Tn}}
   -(-1)^{\alpha\cdot\mu+\alpha\cdot\nu+\mu\cdot\nu}\cdot
    \pstree[treemode=R,levelsep=0.8cm, treesep=0.3cm]{\Tr*{1}}
    {\pstree{\Tc{3pt}} {\pstree{\Tc*{3pt}}
    { \Tcircle{1} \Tcircle{2} \Tp[edge={\ncline[linewidth=0.1cm]}]}}}
\end{multline*}
A thorough investigation shows that all terms in $\delta(H')$ cancel, except for
\begin{multline*}
 \left( (-1)\cdot\pstree[treemode=R, levelsep=0.8cm, treesep=0.3cm,
     linewidth=0.9pt] {\Tcircle{2}}{\pstree[levelsep=0cm]{\Tc{3pt}}
  { \pstree[treemode=U, levelsep=0.5cm]{\Tn} {\Tcircle{1}}
    \pstree[treemode=R, levelsep=0.8cm]{\Tn} {\Tp[edge={\ncline[linewidth=0.1cm]}]}
    \Tn}}-(-1)^{\mu\cdot\nu}\cdot
  \pstree[treemode=R, levelsep=0.8cm, treesep=0.3cm,
     linewidth=0.9pt] {\Tcircle{2}}{
    \pstree[treemode=R, levelsep=0.8cm]{\Tc{3pt}}
       {\pstree{\Tcircle{1}}
       {\Tp[edge={\ncline[linewidth=0.1cm]}]}}}
    -(-1)^{\mu\cdot\nu}\cdot
\pstree[treemode=R, levelsep=0.8cm, treesep=0.3cm,
     linewidth=0.9pt] {\Tcircle{2}}{\pstree[levelsep=0cm]{\Tc{3pt}}
  { \Tn
    \pstree[treemode=R, levelsep=0.8cm]{\Tn} {\Tp[edge={\ncline[linewidth=0.1cm]}]}
    \pstree[treemode=D, levelsep=0.5cm]{\Tn} {\Tcircle{1}}}}\right) \\
  -(-1)^{\mu\cdot\nu}\cdot \left(
 (-1)\cdot\pstree[treemode=R, levelsep=0.8cm, treesep=0.3cm,
     linewidth=0.9pt] {\Tcircle{1}}{\pstree[levelsep=0cm]{\Tc{3pt}}
  { \pstree[treemode=U, levelsep=0.5cm]{\Tn} {\Tcircle{2}}
    \pstree[treemode=R, levelsep=0.8cm]{\Tn} {\Tp[edge={\ncline[linewidth=0.1cm]}]}
    \Tn}}-(-1)^{\mu\cdot\nu}\cdot
  \pstree[treemode=R, levelsep=0.8cm, treesep=0.3cm,
     linewidth=0.9pt] {\Tcircle{1}}{
    \pstree[treemode=R, levelsep=0.8cm]{\Tc{3pt}}
       {\pstree{\Tcircle{2}}
       {\Tp[edge={\ncline[linewidth=0.1cm]}]}}}
    -(-1)^{\mu\cdot\nu}\cdot
\pstree[treemode=R, levelsep=0.8cm, treesep=0.3cm,
     linewidth=0.9pt] {\Tcircle{1}}{\pstree[levelsep=0cm]{\Tc{3pt}}
  { \Tn
    \pstree[treemode=R, levelsep=0.8cm]{\Tn} {\Tp[edge={\ncline[linewidth=0.1cm]}]}
    \pstree[treemode=D, levelsep=0.5cm]{\Tn} {\Tcircle{2}}}}
    \right),
\end{multline*}
which is equal to $Z(f,g)-(-1)^{\mu\cdot\nu}\cdot Z(g,f)$ when applied to elements $f, g\in C^\bullet(A,A)$, together with the remaining terms
\begin{multline*}
  (-1)^{\mu+\nu+\mu\cdot\nu}\cdot
   \pstree[treemode=R, levelsep=0.8cm, treesep=0.3cm,
     linewidth=0.9pt] {\Tr*{1}}{\pstree[levelsep=0cm]{\Tc{3pt}}
  { \pstree[treemode=U, levelsep=0.4cm]{\Tn} {
    \pstree[levelsep=0.7cm]{\Tc*{3pt}}{\Tcircle{1}
    \Tcircle{2}}}
    \pstree[treemode=R, levelsep=0.8cm]{\Tn}
    {\Tp[edge={\ncline[linewidth=0.1cm]}]}\Tn}}
    -(-1)^{\mu\cdot\nu} \cdot \left(
   (-1)^{\alpha\cdot\mu+\alpha\cdot\nu}\cdot
    \pstree[treemode=R,levelsep=0.8cm, treesep=0.3cm]{\Tr*{1}}
    {\pstree{\Tc{3pt}} {\pstree{\Tc*{3pt}}
    { \Tcircle{1} \Tcircle{2}
    \Tp[edge={\ncline[linewidth=0.1cm]}]}}} \right. \\ +
    \pstree[treemode=R,levelsep=0.8cm, treesep=0.3cm]{\Tr*{1}}
    {\pstree{\Tc{3pt}} {\pstree{\Tc*{3pt}}
    { \Tp[edge={\ncline[linewidth=0.1cm]}] \Tcircle{1}
    \Tcircle{2}}}}
   + (-1)^{\alpha\cdot\mu}\cdot
   \pstree[treemode=R, levelsep=0.8cm, treesep=0.3cm,
     linewidth=0.9pt] {\Tr*{1}}{\pstree[treemode=R, levelsep=0.8cm]{\Tc{3pt}}
       {\pstree{\Tc*{3pt}} {\Tcircle{1}
       \pstree{\Tcircle{2}}{\Tp[edge={\ncline[linewidth=0.1cm]}]}
       }}}  \\ \left.
  + (-1)^{\alpha\cdot\mu}\cdot
    \pstree[treemode=R,levelsep=0.8cm, treesep=0.3cm]{\Tr*{1}}
    {\pstree{\Tc{3pt}} {\pstree{\Tc*{3pt}}
    {\Tcircle{1} \Tp[edge={\ncline[linewidth=0.1cm]}]
    \Tcircle{2}}}}
 +\pstree[treemode=R, levelsep=0.8cm, treesep=0.3cm,
     linewidth=0.9pt] {\Tr*{1}}{\pstree[treemode=R, levelsep=0.8cm]{\Tc{3pt}}
       {\pstree{\Tc*{3pt}} {\pstree{\Tcircle{1}}{\Tp[edge={\ncline[linewidth=0.1cm]}]}
         \Tcircle{2}}}} \right) -(-1)^{\mu\cdot\nu}\cdot
\pstree[treemode=R, levelsep=0.8cm, treesep=0.3cm,
     linewidth=0.9pt] {\Tr*{1}}{\pstree[levelsep=0cm]{\Tc{3pt}}
  { \Tn
    \pstree[treemode=R, levelsep=0.8cm]{\Tn} {\Tp[edge={\ncline[linewidth=0.1cm]}]}
    \pstree[treemode=D, levelsep=0.4cm]{\Tn} {
    \pstree[levelsep=0.7cm]{\Tc*{3pt}}{\Tcircle{2}
    \Tcircle{1}}}}},
\end{multline*}
which are equal to 
\begin{multline*}
-(-1)^{\mu\cdot\nu}\cdot \left( (-1)^{\mu+\nu+1}\cdot
    \pstree[treemode=R, levelsep=1cm, treesep=0.3cm]{\Tr*{1}}
    { \pstree[levelsep=0cm]{\Tc{3pt}}
  { \pstree[treemode=U, levelsep=0.5cm]{\Tn}
     {\Tr*{M(g,f)}}
    \pstree[treemode=R, levelsep=1cm]{\Tn}
     {\Tp[edge={\ncline[linewidth=0.1cm]}]}
    \Tn
  }}+
    \pstree[treemode=R, levelsep=1cm, treesep=0.3cm]
    {\Tr*{1}} {\pstree{\Tc{3pt}}
    {\pstree[linewidth=0.1cm]{\Tr*{M(g,f)}}{\Tp} }}
    +
    \pstree[treemode=R, levelsep=1cm, treesep=0.3cm]{\Tr*{1}}
    { \pstree[levelsep=0cm]{\Tc{3pt}}
  { \Tn
    \pstree[treemode=R, levelsep=1cm]{\Tn}
     {\Tp[edge={\ncline[linewidth=0.1cm]}]}
    \pstree[treemode=D, levelsep=0.5cm]{\Tn}
     {\Tr*{M(g,f)}}  }} \right).
\end{multline*}
In the last line we indicated that for $f,g\in C^\bullet(A,A)$, we apply $M(g,f)\in C^\bullet(A,A)$ to the operations associated to the sum of symbols. Using lemma \ref{2x43}, and the graded commutativity of $M$, we see that this is exactly $ -(-1)^{\mu\cdot\nu}\cdot\Delta(M(g,f))\cong -\Delta(M(f,g))$.
\end{enumerate}
\end{proof}

\end{document}